\documentclass[11pt]{amsart}
\usepackage{pb-diagram,amssymb,epic,eepic,verbatim}
\usepackage{hyperref,graphics,epsfig,psfrag}
\hypersetup{colorlinks}
\usepackage{color}
\definecolor{darkred}{rgb}{0.5,0,0}
\definecolor{darkgreen}{rgb}{0,0.5,0}
\definecolor{darkblue}{rgb}{0,0,0.5}
\hypersetup{ colorlinks,
linkcolor=darkblue,
filecolor=darkgreen,
urlcolor=darkred,
citecolor=darkblue }
\newtheorem{theorem}{Theorem}[section]

\newtheorem{corollary}[theorem]{Corollary}

\newtheorem{proposition}[theorem]{Proposition}
\newtheorem{lemma}[theorem]{Lemma}
\newtheorem{lem}[theorem]{}
\theoremstyle{definition}
\newtheorem{definition}[theorem]{Definition}
\theoremstyle{remark}
\newtheorem{remark}[theorem]{Remark}
\newtheorem{example}[theorem]{Example}
\newcommand{\blem}{\begin{lem} \rm}
\newcommand{\elem}{\end{lem}}

\newcommand\cA{\mathcal{A}}
\newcommand\cS{\mathcal{S}}
\newcommand\cO{\mathcal{O}}
\newcommand\M{\mathcal{M}}
\newcommand\D{\mathcal{D}}

\renewcommand\M{\mathcal{M}}
\renewcommand\D{\mathcal{D}}
\newcommand\MM{\mathfrak{M}}

\newcommand{\R}{\mathbb{R}}

\newcommand{\RR}{\mathcal{R}}
\newcommand{\C}{\mathbb{C}}
\newcommand{\cC}{\mathcal{C}}

\newcommand{\Z}{\mathbb{Z}}
\newcommand{\Q}{\mathbb{Q}}

\renewcommand{\P}{\mathbb{P}}
\newcommand{\HH}{\mathbb{H}}

\newcommand\lie[1]{\mathfrak{#1}}

\newcommand{\g}{\lie{g}}

\newcommand{\q}{\lie{q}}

\renewcommand{\t}{\lie{t}}

\renewcommand{\u}{\lie{u}}

\newcommand{\on}{\operatorname}

\newcommand{\fr}{{\on{fr}}}

\newcommand{\dual}{\vee}

\newcommand{\Fr}{\on{Fr}}

\newcommand{\End}{\on{End}}

\newcommand{\fin}{\on{fin}}

\newcommand{\Aut}{ \on{Aut} }

\newcommand{\Ad}{ \on{Ad} }

\newcommand{\Rep}{\on{Rep}}
\newcommand{\Resid}{\on{Resid}}
\newcommand{\Restr}{\on{r}}
\newcommand{\Hom}{ \on{Hom}}
\newcommand{\Ind}{ \on{Ind}}
\renewcommand{\ker}{ \on{ker}}

\newcommand{\Spec}{  \on{Spec}}

\newcommand\dirac{/\kern-1.2ex\partial} 
\newcommand\qu{/\kern-.7ex/} 
\newcommand\hqu{/\kern-.7ex/\kern-.7ex/\kern-.7ex/}
\newcommand\lqu{\backslash \kern-.7ex \backslash} 

\newcommand\dr{r_+ \kern-.7ex - \kern-.7ex r_-}
 
\newcommand{\lev}{{\on{lev}}} 

\newcommand{\labell}\label

\newcommand{\ol}{\overline}

\newcommand\Phinv{\Phi^{-1}}

\newcommand\eps{\epsilon}

\newcommand{\lan}{\langle}
\newcommand{\ran}{\rangle}

\newcommand{\ti}{\tilde}

\newcommand\pt{\on{pt}}

\newcommand\cZ{\mathcal{Z}}

\renewcommand{\ss}{{\on{ss}}}

\newcommand\Gr{\on{Gr}}
\newcommand\Map{\on{Map}}
\newcommand\rank{\on{rank}}
\newcommand\ev{\on{ev}}
\newcommand\Eul{\on{Eul}}

\newcommand\quot{\on{quot}}

\newcommand\ul{\underline}
\newcommand\mO{\mathcal{O}}

\newcommand\bra[1]{ < \kern-.7ex {#1} \kern-.7ex >} 
\newcommand\bdefn{\begin{definition}}
\newcommand\edefn{\end{definition}}
\newcommand\bea{\begin{eqnarray*}}
\newcommand\eea{\end{eqnarray*}}
\newcommand\bcv{\left[ \begin{array}{r} }
\newcommand\ecv{\end{array} \right] }

\newcommand\bma{\left[ \begin{array}{l} }
\newcommand\ema{\end{array} \right]}
\newcommand\ben{\begin{enumerate}}
\newcommand\een{\end{enumerate}}
\newcommand\beq{\begin{equation}}
\newcommand\eeq{\end{equation}}
\newcommand\bex{\begin{example}}
\newcommand\bsj{\left\{ \begin{array}{rrr} }
\newcommand\esj{\end{array} \right\}}

\newcommand\Ch{\on{Ch}}

\newcommand\eex{\end{example}}

\newcommand\sx{*\kern-.5ex_X}

\newcommand{\Hilb}{\on{Hilb}}

\def\mathunderaccent#1{\let\theaccent#1\mathpalette\putaccentunder}
\def\putaccentunder#1#2{\oalign{$#1#2$\crcr\hidewidth \vbox
to.2ex{\hbox{$#1\theaccent{}$}\vss}\hidewidth}}

\begin{document}

\title[Quantum Witten localization]{Quantum Witten localization and \\ 
abelianization for qde solutions}

\author{Eduardo Gonz\'alez} 

\address{
Department of Mathematics
University of Massachusetts Boston
100 William T. Morrissey Boulevard
Boston, MA 02125}
  \email{eduardo@math.umb.edu}

\author{Chris T. Woodward}

\address{Mathematics-Hill Center,
Rutgers University, 110 Frelinghuysen Road, Piscataway, NJ 08854-8019,
U.S.A.}  \email{ctw@math.rutgers.edu}

\thanks{Partially supported by grants DMS1104670 and DMS1207194.  A
  previous version was titled {\em Area-dependence in gauged
    Gromov-Witten theory.}  }

\begin{abstract}  
We prove a quantum version of the localization formula of Witten
\cite{wi:tw}, see also \cite{te:qu}, \cite{pa:mo}, \cite{wood:norm},
that relates invariants of a git quotient with the equivariant
invariants of the action.  Using the formula we prove a quantum
version of an abelianization formula of S. Martin \cite{mar:sy},
relating invariants of geometric invariant theory quotients by a group
and its maximal torus, conjectured by Bertram, Ciocan-Fontanine, and
Kim \cite{be:qu}.  By similar techniques we prove a quantum
Lefschetz principle for holomorphic symplectic reductions.  As an
application, we give a formula for the fundamental solution to the
quantum differential equation (qde) for the moduli space of points on
the projective line and for the smoothed moduli space of framed
sheaves on the projective plane (a Nakajima quiver variety).
\end{abstract}
\maketitle

\parskip 0in
\tableofcontents
\parskip .1in

\section{Introduction}

\subsection{Quantum Witten localization} 

The main result of this paper is a formula relating the equivariant
Gromov-Witten graph invariants of a smooth projective variety with
group action and the graph invariants of the geometric invariant
theory quotient. As a consequence we obtain versions of ``quantum
abelianization'' for graph Gromov-Witten invariants as well as a
``quantum Lefschetz'' principle for holomorphic symplectic reductions.

To state the main result we introduce the following notation. Let $G$
be a connected complex reductive group acting on a smooth polarized
projective variety $X$.  Let $X \qu G$ denote the {\em git quotient}
of $X$ by $G$, which here means the stack-theoretic quotient of the
semistable locus by the group action.  We assume that $G$ acts with
only finite stabilizers on the semistable locus.  In this case the git
quotient $X \qu G$ is a smooth proper Deligne-Mumford stack with
projective coarse moduli space by Mumford et al \cite{mu:ge}.  Let
$H(X \qu G)$ resp. $H_G(X)$ denote the rational resp. equivariant
rational cohomology of $X \qu G$ resp. $X$.  Kirwan's thesis
\cite{ki:coh} studies the natural map
$$ \kappa_{X,G} : H_G(X) \to H(X \qu G) $$
given by restriction to the semistable locus and descent.
Integration over $X \qu G$ defines a {\em trace map}
$$ \tau_{X \qu G}: H(X \qu G) \to \Q, \quad h \mapsto
\int_{[ X \qu G]} h .$$

Naturally one wants to compute the composition of the trace with
Kirwan's surjection.  For example, one would like to compute the
cohomology of $X \qu G$ in terms of the $G$-equivariant cohomology of
$X$.  Witten \cite{wi:tw} introduced a strategy, which he termed {\em
  non-abelian localization}, to compute the composition
$\tau_{X \qu G} \circ \kappa_{X,G}$.  His formula involves a trace map
$$ \tau_X^G: H_G(X) \to \Q, \quad h \mapsto \int_{X \times \g_\R}
h
$$
 given by integration over $X$ {\em and} the unitary part $\g_\R$ of
 the Lie algebra using suitable regularization procedures
 \cite{pa:loc}, \cite{pa:mo}, \cite{wood:norm}.  In the $K$-theory
 version discussed in Paradan \cite{pa:lo}, the Witten trace is easier
 to define: it is the invariant part of the index, and no
 regularization procedure is needed.  Witten's localization formula
 computes the difference between $\tau_{X \qu G} \circ \kappa_{X,G}$
 and $\tau_X^G$, that is, the failure of the following diagram to
 commute:
\begin{equation} \label{naldiag} \begin{diagram} \node{H_G(X)} \arrow{se,b}{\tau_X^G} \arrow[2]{e,t}{\kappa_{X,G}}
\node{} \node{H(X \qu G)} \arrow{sw,r}{\tau_{X \qu G}} \\ \node{}
\node{\Q} \node{} \end{diagram} .\end{equation}
By Witten's argument in \cite{wi:tw}, \cite{pa:loc} the difference is
a sum of contributions from the Kirwan-Ness strata of positive
dimension:
$$ \tau_X^G = \tau_{X \qu G} \circ \kappa_{X,G} + \sum_{ [\zeta] \neq
  0}  \tau_{X,G, \zeta}: H_G(X) \to \Q $$
where $\tau_{X,G,\zeta}$ is a contribution from the stratum with
maximally-destabilizing one-parameter subgroup $\zeta$.  An explicit
formula for the contributions $\tau_{X,G,\zeta}$ was described in
papers by Teleman \cite{te:qu} in the case of sheaf cohomology, by
Paradan \cite{pa:lo} for $K$-theory of Hamiltonian actions, and in
papers by Paradan and Woodward \cite{pa:loc}, \cite{pa:mo},
\cite{wood:norm} for cohomology of Hamiltonian actions, see also
Beasley-Witten \cite{bw:nal} which uses the localization formula to
compute the Chern-Simons partition function for Seifert manifolds.  A
different formula computing the composition is given in Jeffrey-Kirwan
\cite{je:lo1}.  A virtual Witten localization formula has recently
appeared in Halpern-Leistner \cite[(5)]{hl:qs}.

The quantum version of Witten's localization formula compares
Gromov-Witten invariants of a git quotient with equivariant
Gromov-Witten invariants for the action.  To state the result let
$\omega \in H_2^G(X)$ be the first Chern class of the linearization
(that is, the symplectic class) and let
$$ \Lambda_X^G = \left\{ \sum_{i =0}^\infty c_i q^{d_i} ,
c_i \in \Q, d_i \in H_2^G(X,\Q), \ \lim_{i \to \infty} \lan
d_i,\omega \ran = \infty \right\} $$
denote the equivariant Novikov field for $X$.  Let
$$ QH_G(X) = H_G(X) \otimes \Lambda_X^G $$
denote the equivariant quantum cohomology of $X$.  Virtual integration
over the moduli stack of $n$-marked genus $0$ stable maps
$\ol{\M}_{0,n}(X)$ for $n \ge 3$ defines a family of formal quantum
products
$$ \star_h: T_h QH_G(X)^2 \to T_h QH_G(X),
\quad h \in QH_G(X) .$$
Formal in this setting means that only the Taylor coefficients of the
maps are convergent.  Define a quantum version of Witten's trace as
follows.  Let $\P = (\C^2 - \{ 0 \})/\C^\times$ denote the projective
line.  For $d \in H_2(X,\Z)$ let
$ \ol{\M}_n(\P,X,d) := \ol{\M}_{0,n}(\P \times X,(1,d)) $
denote the moduli stack of parametrized stable maps from $\P$ to $X$
of class $d \in H_2^G(X,\Z)$.  The action of $G$ on $X$ induces a
natural action on $\ol{\M}_n(\P,X,d)$.  A natural stability condition
for the action is given by requiring that the stable map has
generically semistable value \cite{small}.  Denote by
$ \ol{\M}_n(\P,X,d) \qu G$ the stack-theoretic quotient of the
semistable locus by the group action.  By, for example, \cite[Lemma
2.6]{wall}, $\ol{\M}_n(\P,X,d) \qu G$ is a proper Deligne-Mumford
stack with a perfect relative obstruction theory.  Via equivariant
formality we may consider $H_2(X,\Z)/\on{torsion}$ as a subgroup of
$H_2^G(X,\Q)$.  Denote by $\tau_X^G$ the formal trace map given by
virtual integration over the moduli stacks $\ol{\M}_n(\P,X,d) \qu G$:
$$ \tau_X^G: QH_G(X) \to \Lambda_X^G, \quad h \mapsto \sum_{n \ge 0, d
  \in H_2(X,\Z)/\on{torsion}} (q^d/n!)  \int_{[\ol{\M}_n(\P,X,d) \qu
    G]} \ev^* (h \otimes \ldots \otimes h) $$
for $h \in H_G(X)$.   The map $\tau_X^G$ is a quantum version of
Witten's trace in the sense that if one sets $q = 0$ and fixes the
positions of the markings then one obtains the classical Witten trace
for polynomial classes, that is, the integral over $X \qu G$.

A quantum version of Kirwan's map counting maps to the quotient stack
with semistability enforced at a marked point was introduced in
\cite{qk1}, \cite{qk2}, \cite{qk3}.  The quantum Kirwan map is a
non-linear map, still denoted $\kappa_{X,G}$,
\begin{equation} \label{qkirwan} \kappa_{X,G}: QH_G(X) \to
  QH(X \qu G) \end{equation}
with the property that any linearization 
$$D_h \kappa_{X,G}: T_h QH_G(X) \to
T_{\kappa_{X,G}(h)} QH(X \qu G)$$ 
is a homomorphism with respect to the quantum products.  In
particular, if $\kappa_{X,G}(0) = 0$ (which generally happens only in
Fano cases) then $D_0 \kappa_{X,G}$ is a homomorphism from the small
equivariant quantum cohomology $T_0 QH_G(X)$ of $X$ to the quantum
cohomology $T_0 QH(X \qu G)$ of $X \qu G$.

A quantum version of the integration over the geometric invariant
theory quotient is defined by a count of stable maps to the graph
space.  Recall that $\ol{\M}_n(\P,X \qu G,d)$ denotes stable maps to
$\P \times (X \qu G)$ of class $1,d$.  Using the Behrend-Fantechi
virtual fundamental classes define
$$ \tau_{X \qu G}: QH(X \qu G) \to \Lambda_X^G, \quad h \mapsto
\sum_{n \ge 0, d \in H_2(X \qu G,\Q)} (q^d/n!)  \int_{[
\ol{\M}_n(\P,X
    \qu G,d)]}
 \ev^* (h \otimes \ldots \otimes h) $$
for $h \in H(X \qu G)$.
The quantum Witten localization formula gives a precise description of
the difference between the traces $\tau_X^G$ and $\tau_{X \qu G} \circ
\kappa_{X,G}$.  That is, it measures the failure of the ``quantum
integration'' to commute with reduction, i.e.  the failure of
commutativity of the diagram
\begin{equation} \label{naldiag2} \begin{diagram}
	\node{QH_G(X)} \arrow{se,b}{\tau_X^G}
	\arrow[2]{e,t}{\kappa_{X,G}} \node{} \node{QH(X \qu G)}
	\arrow{sw,r}{\tau_{X \qu G}} \\ \node{}
	\node{\Lambda_X^G} \node{} \end{diagram} .\end{equation}
    As in the classical Witten localization formula \cite{wi:tw}, the
    failure to commute is given by a sum of fixed point contributions.
    Each term is a {\em gauged Gromov-Witten invariant}
    $\tau_{X,G,\zeta,\rho}$ associated to the action of centralizers
    on components of the fixed point variety of some one-parameter
    subgroup $\exp(\C \zeta) \subset G$, stable with respect to the
    linearization $\ti{X}^\rho$ for some $\rho \in (0,\infty)$.  The
    main result is the following:

    \begin{theorem} \label{qwform} {\rm (Quantum Witten localization)}
      Let $C$ be a smooth connected projective curve of genus $0$, $X$
      a smooth projective $G$-variety, and $\ti{X}$ a linearization.
      Suppose that for every $\zeta \in \g$ and $\rho \in (0,\infty)$,
      stable=semistable for $\ol{\M}_n^{G}(C,X,\ti{X}^\rho,\zeta)$,
      and stable=semistable for the $G$-action on $X$.  Then the
      following equality holds for formal maps from $QH_G(X)$ to
      $ \Lambda_X^G$:
\begin{equation} \label{qwitten}   \tau_X^G - \tau_{X \qu G}
  \circ \kappa_{X,G} = \sum_{ [\zeta] \neq 0, \rho \in (0,\infty) }
 \tau_{X,G,\zeta,\rho} 
  .\end{equation}
\end{theorem} 

\subsection{Applications to quantum abelianization}

We give two groups of applications.  The first group consist of
versions of the {\em quantum Martin conjecture} of Bertram et al
\cite{be:qu} that compares Gromov-Witten invariants of a git quotient
$X \qu G$ and the quotient $X \qu T = X^{\ss,T}/T$ by a maximal torus
$ T \subset G$.  As an example, we give a formula for the qde solution
for the quotient of points on the projective line by its
automorphisms, see \eqref{first} \eqref{second} \eqref{third} below.

Before describing the quantum generalization we discuss the classical
story of abelianization due to Martin \cite{mar:sy}. Let $\nu_{\g/\t}$
denote the bundle over $X \qu T$ induced from the trivial bundle with
fiber $\g/\t$ over $X$ and $\tau_{X \qu T}^{\g/\t} $ the
Euler-twisted integration map
$$ \nu_{\g/\t} = X^{\ss,T} \times_T (\g/\t), \quad \tau_{X \qu T}^{\g/\t} 
: H(X \qu T) \to \Q, \quad h \mapsto \int_{ [X \qu T]} h \cup
\Eul(\nu_{\g/\t}) .$$
Let $W = N(T) /T$ denote the Weyl group of $T \subset G$ and
$\Restr_T^G$ the isomorphism with Weyl-invariants
$$\Restr_T^G: H_G(X) \cong H_T(X)^W.$$

\begin{theorem} \label{martin}  {\rm (Martin formula
	\cite{mar:sy})}    Let $X$ be a smooth projective
	$G$-variety.  Suppose that stable=semistable for the
	actions of $T$ and $G$ on $X$.  Then integration over $X
	\qu G$ and $X \qu T$ are related by
$$ \tau_{X \qu G} \circ \kappa_{X,G} = |W|^{-1} \tau^{\g/\t}_{X \qu
T} \circ \kappa_{X,T} \circ \Restr_T^G 
.$$
Furthermore, there exists a surjective map
\begin{equation} \label{muTG} \mu_T^G : H(X \qu T)^{W} \to H(X
  \qu G) \end{equation}
whose kernel is the annihilator of $\Eul(\nu_{\g/\t})$. 
\end{theorem}  
\noindent The first part of the theorem is equivalent to commutativity
of the following diagram:
$$ \hskip .65in \begin{diagram} \node{} \node{H_G(X)
  \cong H_T(X)^W} \arrow{se,l}{\kappa_{X,T}}
  \arrow{sw,l}{\kappa_{X,G}} \node{} \\ \node{H(X \qu G)}
  \arrow{se,r}{\tau_{X \qu G}} \node{} \node{H(X \qu T)}
  \arrow{sw,r}{ |W|^{-1} \tau_{X \qu T}^{\g/\t}} \node{} \\
  \node{} \node{\Q} \end{diagram} $$

Using quantum Witten localization \eqref{qwitten} we prove
that a formula similar to that in Theorem \ref{martin} holds in
quantum cohomology.  Versions of this formula were conjectured and
several special cases proved by Hori-Vafa \cite[Appendix]{ho:mi},
Bertram-Ciocan-Fontanine-Kim \cite{be:qu} and
Ciocan-Fontanine-Kim-Sabbah \cite{ciocan:abnonab}.  A general result
that holds under monotonicity conditions is proved by Schm\"ashcke
\cite{schm:thes}.  The push-forward in homology $ \pi_T^G : H_2^T(X)
\to H_2^G(X)$ defines a map of equivariant Novikov rings
\begin{equation} \label{piTG} \pi_T^G : \Lambda_X^T \to
  \Lambda_X^G, \quad \sum_{d \in H_2^T(X)} c_d q^d \mapsto
  \sum_{d \in H_2^G(X)} c_d q^{\pi(d)} .\end{equation}
Let $Q\HH_G(X) \subset QH_G(X)$ denote the
subspace generated by Chern characters of algebraic vector
bundles,
$$ \HH_G(X) := \{ \Ch_G(E) \ | \ E \to X \text{ \ vector bundle } \},
\quad Q\HH_G(X) := \HH_G(X) \otimes \Lambda_X^G .$$
The restriction to Chern characters is necessary because our arguments
at some point use sheaf cohomology.  We denote by
$$\Restr_T^G: QH_G(X) \to QH_T(X)$$ 
the map obtained by combining the pull-back $H_G(X) \to H_T(X)$ with
the inclusion $\Lambda_X^G \subset \Lambda_X^G$ induced by the
inclusion $H_2^G(X,\Z) \cong H_2^T(X,\Z)^W \subset H_2^T(X,\Z)$.

\begin{theorem} {\rm (Quantum Martin formula)} \label{qabel}
  Let $C$ be a smooth connected projective genus $0$ curve and $X$ a
  smooth linearized projective $G$-variety.  Suppose that
  stable=semistable for $T$ and $G$ actions on $X$.  The following
  equality holds on $Q\HH_G(X)$:
\begin{eqnarray*}
 \tau_{X \qu G} \circ \kappa_{X,G} &=& |W|^{-1} \pi_T^G \circ 
\tau_{X
    \qu T}^{\g/\t} \circ \kappa_{X,T}^{\g/\t} \circ \Restr_T^G \\
&=&
  |W|^{-1} \pi_T^G \circ \tau_{X, T}^{\g/\t} \circ \Restr_T^G:
  \ \ \ Q\HH_G(X) \to \Lambda_X^G. \end{eqnarray*}
\end{theorem} 

\noindent That is, there is a commutative diagram
$$ \hskip .6in \begin{diagram}
  \node{Q\HH_G(X)}\arrow{s,l}{\kappa_{X,G}}
  \arrow{e} \node{QH_T(X)} \arrow{s,r}{\kappa_{X,T}^{\g/\t}} \\
  \node{QH(X \qu G)} \arrow{s,l}{\tau_{X \qu G}} \node{QH(X
  \qu T)} \arrow{s,r}{|W|^{-1} \tau_{X \qu T}^{\g/\t}} \node{} \\
  \node{\Lambda_X^G} \node{\Lambda_X^T}  \arrow{w}
\end{diagram}$$

A similar abelianization formula holds for solutions to quantum
differential equations.  The $\C^\times$-equivariant extension of the
graph potential on a genus zero curve (with the standard
$\C^\times$-action by rotations) admits a factorization into {\em
  localized graph Gromov-Witten potentials}
$$ \tau_{X \qu G,\pm}: QH(X \qu G) \to QH_{\C^\times}(X \qu G) .$$
Here, as in the remainder of the paper, $QH_{\C^\times}$ denotes the
completion of the $\C^\times$-equivariant cohomology with the
equivariant generator inverted and completed: If the equivariant
generator is denoted $\zeta \in H(B\C^\times)$ then
$$ QH_{\C^\times}(X \qu G) = QH(X \qu G)[\zeta,\zeta^{-1}]] .$$ 
In the literature these potentials are often call $J$-functions or
one-point descendant potentials \cite{gi:eq}.  Denote by
$$ \mu_T^G : QH(X \qu G) \to QH(X \qu T) $$
the map combining Martin's map of \eqref{muTG} with the canonical map
of Novikov rings $\pi_T^G$ of \eqref{piTG}.

\begin{theorem} \label{abelloc}
 {\rm (Abelianization for qde solutions)} Suppose that $X$ is a smooth
 linearized projective $G$-variety, and stable=semistable for the $T$
 and $G$-actions on $X$.  Then
\begin{eqnarray*} 
 \tau_{X \qu G,\pm} \circ \kappa_{X,G} &=& \mu_T^G \circ\tau^{\g/\t}_{X
   \qu T, \pm} \circ \kappa_{X,T}^{\g/\t} \circ \Restr_T^G \\
&=& \mu_T^G
 \circ\tau^{\g/\t}_{X,T,\pm} \circ \Restr_T^G: Q\HH_G(X) \to QH_{\C^\times}(X \qu
 G) .\end{eqnarray*} 
 \end{theorem}

The argument extends to quasiprojective targets under suitable
properness conditions.  In particular, it holds for targets that are
$G$-vector spaces $X$ satisfying a certain convexity condition, see
Theorem \ref{vspace} below.  In Examples \ref{grass}, \ref{points} we
apply the formula to give formulas for the solution to the quantum
differential equation for the Grassmannians and moduli of points on
the projective line.

\begin{example} \label{points} {\rm (Moduli of points on the projective line) } 
We consider the git quotient for the diagonal action of $SL(2,\C)$ on
$(\P^1)^{2k+1}$ with linearization on each factor the same.  The git
quotient is
$$ Y = \left\{ (x_1,\ldots, x_{2k+1}) \in (\P^1)^{2k+1} \ | \ \sup_{x \in
  \P} \# \{ x_i = x \} \leq k \right\} / SL(2,\C) .$$
%
In order to apply our results we realize the moduli space of points as
a git quotient of a vector space.  The product $(\P^1)^{2k+1}$ is the
git quotient of $X = \C^{4k+2}$ by the diagonal action of
$(\C^\times)^{2k+1}$.  Thus,
$$ Y = X \qu G, \quad X = \C^{4k+2}, \quad G = (\C^\times)^{2k+1}
\times SL(2,\C) .$$
The diagonal subgroup $\C^\times$ in the first factor acts on $X$
with all positive weights, so $X$ is convex.  Since $X$ is
$G$-equivariantly Fano the quantum Kirwan map is the identity on
$QH^{\leq 2}_G(X)$ for reasons of dimension by \eqref{eqfano}.  By
abelianization the localized graph potential is given by
\begin{equation} \label{first} \tau_{X \qu G,\pm} =
\mu_T^G \circ \tau_{X,T,\pm} \circ \Restr_T^G \in QH_{\C^\times} (X \qu
G) .\end{equation}
The maximal torus of $G$ is $T \cong (\C^\times)^{2k+2}$ embedded as
the subgroup of products of diagonal matrices.  The canonical
identifications as symmetric polynomials gives
$$ H_2^G(X,\Z) \cong H_2^{(\C^\times)^{2k+1}}(X,\Z) \cong \Z^{2k+1}, \quad
H_2^T(X,\Z) \cong H_2^{(\C^\times)^{2k+2} }(X,\Z) \cong \t_\Z \cong
\Z^{2k+2} .$$
The weights for the $T$-action on $X$ are written in terms of the
standard basis $\eps_1,\ldots,\eps_{2k+2}$
$$ \eps_1 + \eps_{2k+2}, \eps_1 - \eps_{2k+2}, \eps_2 + \eps_{2k+2},
\eps_2 - \eps_{2k+2} , \ldots \eps_{2k+1} - \eps_{2k+2} \in
\t_\Z^\dual.$$
Let $\theta_1,\ldots, \theta_{2k+2} \in H^2_T(X)$ denote the
generators corresponding to the splitting $T = (\C^\times)^{2k+2}$.
For any $2k+2$-tuple of non-negative integers $\ul{d} =
(d_1,\ldots,d_{2k+2})$, $\theta = \sum c_i \theta_i$ with $c_i \in
\Z$, define a factorial-like product 
$$ \Delta_{\ul{d}}(\theta) := \frac{ \prod_{l = -\infty}^{\theta \cdot \ul{d}} (\theta
  + l \zeta)} { \prod_{l = -\infty}^{0} (\theta + l \zeta)} .$$
The localized potential $\tau_{X,T,,\pm}$ has restriction 
$$ \tau_{X,T,\pm} | QH^{\leq 2}(X \qu G) \subset QH^{\leq 2}_G(X)
\subset QH^{\leq 2}_T(X)$$ 
given by
\begin{multline} \label{second} (\tau_{X,T,\pm} | QH^{\leq 2}(X \qu G)):  QH^{\leq 2}(X \qu G) \to 
QH^{\leq 2}(X \qu G)[[\zeta^{-1}]] \\ \quad (t_0 +t_1 \theta_1 +
\ldots + t_{2k+2} \theta_{2k+2}) \mapsto \sum_{\ul{d}} q^{\ul{d}}
e^{t_0 + (t_1 (\theta_1 + d_1 \zeta) + \ldots + t_{2k+2} (
  \theta_{2k+2} + d_{2k+2}\zeta)) / \zeta}
\tau_{X,T,\pm}(\ul{d}) \end{multline}
where 
\begin{multline} \label{third}
 \tau_{X,T,\pm}(\ul{d}) := e^{\sum_{j=1}^{2k+2} d_j t_j }
 \frac{\Delta_{\ul{d}}( 2 \theta_{2k+2} ) \Delta_{\ul{d}}( - 2
   \theta_{2k+2} ) } {\Delta_{\ul{d}}( \theta_1 + \theta_{2k+2})
   \Delta_{\ul{d}}( \theta_1 - \theta_{2k+2}) \ldots \Delta_{\ul{d}}(
   \theta_{2k+1} -\theta_{2k+2}) } .\end{multline}
Combining \eqref{first}, \eqref{second}, \eqref{third} gives a formula
for a qde solution.
\end{example} 

Using the result on abelianization of qde solutions we obtain the
following relationship between quantum cohomology rings of abelian and
non-abelian quotients which was proved in monotone cases by
Schm\"aschke \cite{schm:thes}.  Let $T \subset G$ a maximal torus as
above.  Denote by 
$$\RR = \RR_+ \cup \RR_- \subset \t^\dual_\Z $$ 
the set $\RR$ of roots, partitioned into positive $\RR_+$ and negative
$\RR_-$ roots.  Consider the decomposition of the Lie algebra into
root spaces,
\begin{equation} \label{rdec} \g \cong \t \oplus \bigoplus_{\alpha \in \RR_-} \g_\alpha \oplus
\bigoplus_{\alpha \in \RR_+} \g_\alpha .\end{equation}
Denote the Euler class of $\g/\t$
$$e = e_- e_+, \quad e_\pm = \prod_{\alpha \in \RR_\pm} \alpha \in
S(\t^\dual) \subset QH_T(X) .$$
The image of $e_\pm$ under the linearized quantum Kirwan map is
denoted
$$ D_h \kappa_{X,T}^{\g/\t}(e_\pm) \in QH(X \qu T), \quad D_h
\kappa_{X,T}^{\g/\t}(e_+) = (-1)^r D_h \kappa_{X,T}^{\g/\t}(e_-).$$
Either $e_-$ or $e_+$ works equally well in the formulas below.

\begin{theorem} \label{compare} 
{\rm (Comparison of quantum cohomology rings)} For any ${h} \in
QH_G(X)$ there exists a canonical surjection
$$ T_{\kappa_{X,T}({h})} QH(X \qu T)^{W} \to
T_{\kappa_{X,G}({h})} QH(X \qu G) $$
whose kernel is the annihilator of $D_h \kappa_{X,T}(e_\pm):$
$$ T_{\kappa_{X,G}(h)} QH(X \qu G) = T_{\kappa_{X,T}({h})} QH(X \qu
T)^{W} / \on{ann} ( D_h \kappa_{X,T}(e_\pm)) .$$
\end{theorem} 

\begin{example} {\rm (Quantum cohomology of the Grassmannian)}  The Grassmannian $\Gr(k,n)$ is the git quotient
$X \qu G$ of $X = \Hom(\C^k,\C^n)$ by the action of $G = GL(k)$.  The
  maximal torus is $T = GL(1)^k$ and the abelian quotient
$$X \qu T \cong (\P^{n-1})^k .$$
The standard presentation of $QH( (\P^{n-1})^k)$ is
$$ QH( (\P^{n-1})^k) = \Lambda_X^G[H_1,\ldots, H_k], \quad H_i^n = q,
i = 1,\ldots, k $$
where $H_i$ is the hyperplane class on the $i$-th factor.  The Weyl
group $W = S_k$ is the $k$-symmetric group.  The $W$-invariant part of
the cohomology ring $QH( (\P^{n-1})^k)$ is generated by the Schur
polynomials
$$ \chi_{\lambda^\dual}(H_1,\ldots,H_n) = \prod_{w \in W} \frac{(-1)^{l(w)}
  H^{w(\lambda + \rho) - \rho}}{ \prod_{i < j} (H_j - H_i) } $$
where 
$$ H^\lambda = H_1^{\lambda_1} \ldots H_k^{\lambda_k}, \quad \rho =
(1,\ldots, k) .$$
Since $X$ is equivariantly Fano $\kappa_{X,T}(0) = 0$ has no quantum
corrections.  The first quantum corrections to $D_0 \kappa_{X,T}$
occur in the twice the degree of the minimal Chern number $n$, so that
$$ D_0 \kappa_{X,T}(e_\pm) = \pm \star_{i < j} (H_j - H_i) .$$
Hence for any $\mu \in \Z^k$, the elements
$$ \chi_{\lambda + n \mu} (H_1,\ldots,H_n) - \chi_{\lambda^\dual}(H_1,\ldots,H_n) $$ 
are in the annihilator of $D_0 \kappa_{X,T}(e_\pm)$.  So we have
relations
$$ \chi_{\lambda + n \mu} = q^{\mu_1 + \ldots + \mu_k} \chi_{\lambda^\dual}
\in QH(\Gr(k,n)), \quad \lambda,\mu \in \Z^k .$$
These are the usual relations in the cohomology of the Grassmannian
describing the cohomology as a truncation of the polynomial
representation ring of $GL(k)$ in for example
Bertram-Ciocan-Fontanine-Fulton \cite{be:qm}.
\end{example} 

\subsection{Applications to holomorphic symplectic quotients} 

A second application of the quantum Witten localization formula is to
the quantum Lefschetz principle, by which Gromov-Witten invariants of
complete intersections are expressed in terms of Euler-twisted
Gromov-Witten invariants of the ambient space.  This extends the
quantum Lefschetz principle beyond cases where the bundle is {\em
  concavex}, that is, a direct sum of convex and concave line bundles
\cite{co:qrr}, \cite{elezi}.  As an example, we give a formula for the
qde solution of the ADHM quiver variety, Theorem \ref{desing} below.

We consider the Gromov-WItten invariants of zero sets of sections of
associated bundles, as follows.  Suppose that $X \qu G$ is a git
quotient as above, and $V$ is a $G$-representation.  Then
$$V \qu G = (V | X^{\ss}) / G \to X \qu G$$ 
is the associated bundle on the git quotient.  Suppose that the bundle
$V \qu G$ admits a section 
$$\Phi: X \qu G \to V \qu G$$ 
induced from an equivariant map $\Phi: X \to V$.  Denote the level
sets
$$Z \qu G = (\Phi \qu G)^{-1}(0) \subset X \qu G, \quad Z :=
\Phinv(0) $$
is a smooth subvariety.  To simplify notation, we denote by
$QH(Z \qu G)$ the quantum cohomology defined over the Novikov ring
$\Lambda_X^G$, that is,
$QH(Z \qu G) = H(Z \qu G) \otimes \Lambda_X^G$.  Let
$r_{Z,G}: QH_G(X) \to QH_G(Z)$ and
$r_{Z \qu G}: QH(X \qu G) \to QH(Z \qu G)$ denote pull-backs.

The following describes the graph potentials of holomorphic symplectic
quotients in terms of the twisted graph potential for the ambient
variety:

\begin{theorem} \label{qlef} 
{(\rm quantum Lefschetz for associated bundles)} Let $X$ be a
linearized projective or convex quasiprojective $G$-variety and $\Phi:
X \to V $ a section as above with smooth zero set $Z$.  The graph
potentials for $Z \qu G$ and $X \qu G$ are related by
$$ \tau_{X \qu G}^{V \qu G} \circ \kappa_{X,G}^V = \tau_{X,G}^V =
 \tau_{Z \qu G} \circ \kappa_{Z,G} \circ r_{Z,G}: QH_G(X) \to \Lambda_X^G .$$
Similarly for the qde solutions
\begin{eqnarray*} 
  r_{X \qu G} \circ \tau_{X \qu G,\pm}^{V \qu G} \circ \kappa_{X,G}^{V}
  &=& r_{X \qu G}^{V \qu G} \circ \tau_{X,G,\pm}^{V}\\ &=& r_{Z \qu G}
                                                           \circ \tau_{Z \qu
                                                           G,\pm} \circ \kappa_{Z,G} \circ r_{Z,G}: QH_G(X) \to QH(X \qu G)
                                                           .\end{eqnarray*}
\end{theorem} 

Combining this result with abelianization allows us to compute the
graph potentials of certain holomorphic symplectic quotients.  Suppose
that $X$ is equipped with an holomorphic moment map
$\Phi: X \to \g^\dual$ as well as a linearization $\ti{X} \to X$.  The
holomorphic symplectic quotient is then the git quotient of the zero
level set:
$$ X \hqu G := \Phinv(0) \qu G = \Phinv(0)^{\ss}/ G .$$
Suppose that $Z = \Phinv(0)$ is smooth.  Let $r^X_Z: Q\HH_G(X) \to
Q\HH_G(Z)$ denote the restriction map.  Let 
$$\mu_{Z \qu G}^{X \qu T}: QH(X \qu T) \to QH(Z \qu G)$$ 
be the combination of pull-back $QH(X \qu G) \to QH(Z \qu G)$ with
Martin's surjection $QH(X \qu T) \to QH(X \qu G)$ \cite{mar:sy}.

\begin{theorem}  \label{hk}  {\rm (qde solutions for holomorphic symplectic quotients)}   
 Let $X$ be a linearized projective or convex quasiprojective
 $G$-variety and $\Phi: X \to \g^\dual $ an equivariant map as above
 with zero set $Z$, and $T \subset G$ a maximal torus.  The graph
 potential for $Z \qu G$ satisfies
$$ \tau_{Z \qu G} \circ \kappa_{Z,G} \circ r_{Z,G} = \tau_{X,T}^{\g
   \oplus \g/\t}: Q\HH_G(X) \to \Lambda_X^G .$$
Similarly the localized graph potentials are related by 
$$ \mu_{Z \qu G}^{X \qu T} \circ \tau_{X,T,\pm}^{\g \oplus \g/\t} =
\tau_{Z \qu G,\pm} \circ \kappa_{Z,G} \circ r_{Z,G}: Q\HH_G(X) \to QH(Z
\qu Q) .$$
\end{theorem}

We apply this formula, at least in principle, to the moduli of framed
sheaves on the projective plane.  Let $\P^2$ denote the projective
plane and let
$$\ell_\infty = \{ [0,z_1,z_2] \} \subset \{ [z_0,z_1,z_2] \} =  \P^2$$ 
denote the divisor at infinity.  Recall
$$ \M = \left\{ (E,\Phi) \left|
\begin{array}{c}
E: \text{torsion free sheaf on $\P^2$} \\
\rank(E) = r, c_2(E) = k \\
\Phi: E|_{\ell_\infty} \to \mO_{\ell_\infty}^{\oplus r}:
\text{framing at infinity} \end{array} \right.  \right\} / \text{isomorphism}  .$$
According to the Atiyah-Drinfeld-Hitchin-Manin description of the
moduli space \cite{adhm}  there exists an isomorphism
\begin{equation} \label{adhmpresent} \M \cong \left\{ 
(B_-,B_+,i_-,i_+) \left|
\begin{array}{c} 
[B_-,B_+] + i_- i_+ = 0 \\
\text{there exists no subspace} \\
S \subset \C^k \text{\ such that\ } B_\pm(S) \subset S \\
\text{\ and\ } \on{im}(i_-) \subset S 
\end{array} \right. \right \} / G ,
\end{equation}
where
$$B_-,B_+ \in \End(\C^k), \quad i_- \in \Hom(\C^r,\C^k), \quad i_+ \in
\Hom(\C^k,\C^r).$$
The action of $g \in G$ is given by
$$ g \cdot (B_-,B_+,i_-,i_+) = (g B_- g^{-1}, g B_+ g^{-1}, g i_- ,
i_+ g^{-1}) .$$
This moduli space is a special case of a Nakajima quiver variety.

The moduli space of gauged maps in the case of the quiver variety
describing the moduli space is not proper.  Instead one introduces an
auxiliary torus action that acts with proper fixed point components
so that the equivariant Gromov-Witten theory is defined by
localization.  The group $S = (\C^\times)^2$ acts equivariantly on
$$ X := \End(\C^k)^{\oplus 2} \oplus \Hom(\C^r,\C^k) \oplus
\Hom(\C^k,\C^r) $$
by 
$$ (s_-,s_+)(B_-,B_+,i_-,i_+) = (s_- B_-, s_+ B_+, i_-, s_- s_+ i_+) .$$
For any $\lambda \in \C$ the action of $S$ preserves the locus
$$ Z = \left\{ (B_-,B_+,i_-,i_+) \left| [B_-,B_+] + i_- i_+ = \lambda \on{Id}
\right. \right\} \subset X $$
and induces an $S$-action on the quotient $\M$.  For any character
$\chi \in \Hom(GL_k(\C),\C^\times) \cong \Z$, let $ \ti{\M}$ denote
the shifted quotient
$$ \ti{\M} = Z \qu_\chi G \subset X \qu_\chi G $$
where $G = GL_k(\C)$ and $\qu_\chi$ denotes the $\chi$-shifted
geometric invariant theory quotient.  We take $\lambda,\chi$ to be
generic small values, so that $\ti{\M}$ is a smooth variety.

To compute the qde solution, let $T \subset G$ denote the diagonal
maximal torus.  We consider the twisted Gromov-Witten theory of $X
\qu_\chi G$ corresponding to the relation defining $Z$, that is,
twisted by the Euler class of the index bundle of $\End(\C^k)$.
Define the factorial-like product for $\theta \in H_2^G(X,\Z) \cong
\Z^k$,
$$ \Delta_{\ul{d}}(\theta,w) := 
\frac{ \prod_{l = -\infty}^{\theta \cdot \ul{d}}  (\theta + w + l \zeta)}
{ \prod_{l = -\infty}^{0}  (\theta + w + l \zeta)} .$$
The twisted localized gauged potential for the $T$ action on
$X$ has restriction to $QH^{\leq 2}_{T}(X)$ given by
(cf. \cite{ciocan:hilb})
$$ \pi_T^G \tau_{X,T,-} = \sum_{d \ge 0} q^d \sum_{\ul{d}: d_1 +
  \ldots + d_k = d} e^{t_0 + (t_1 (\theta_1 + d_1 \zeta) + \ldots +
  t_k \theta_k + d_k \zeta) / \zeta} \tau_{X,T,-}(\ul{d}) $$
where 
$$ \tau_{X,T,-}(\ul{d}) = \prod_{i \neq j}
\frac{\Delta_{\ul{d}}(\theta_i - \theta_j, \xi_- + \xi_+) \Delta_{\ul{d}}(\theta_i
 -\theta_j, 0 )} {\Delta_{\ul{d}}(\theta_i - \theta_j,\xi_-) \Delta_{\ul{d}}(\theta_i -
 \theta_j, \xi_+)} \prod_{i =1 }^k \frac{1}{\Delta_{\ul{d}}(\theta_i,0)^r
 \Delta_{\ul{d}}(-\theta_i, \xi_- + \xi_+)^r }$$
and $\xi_-,\xi_+$ are the equivariant parameters for $S$, that is,
$H_S(\pt) \cong \Q[\xi_-,\xi_+]$.   

\begin{theorem} \label{desing} The localized graph potential of the
  smoothed moduli space of framed sheaves on the projective plane on
  $QH^2_T(X)$ is given by
$$ \tau_{Z \qu G,-} \circ r_{Z,G}
= \mu_{Z \qu G}^{X \qu T} \exp( - \tau_{X,T,-}^{(1)}/\zeta)
\tau_{X,T,-} $$
where $\tau_{X,T,-}^{(1)}$ is the $\zeta^{-1}$-coefficient in 
\begin{multline} \tau_{X,T,-} = \sum_{d \ge 1} q^d \sum_{\ul{d}: d_1 + \ldots + d_k
   = d} e^{t_0 + (t_1 (\theta_1 + d_1 \zeta) + \ldots + t_k (\theta_k
    + d_k \zeta)) / \zeta} \prod_{i \neq j}
  \\ \frac{\Delta_{\ul{d}}(\theta_i - \theta_j, \xi_- + \xi_+)
    \Delta_{\ul{d}}(\theta_i -\theta_j, 0 )} {\Delta_{\ul{d}}(\theta_i
    - \theta_j,\xi_-) \Delta_{\ul{d}}(\theta_i - \theta_j, \xi_+)}
  \prod_{i =1 }^k \Delta_{\ul{d}}(\theta_i,0)^{-r}
  \Delta_{\ul{d}}(-\theta_i, \xi_- + \xi_+)^{-r} .\end{multline}
\end{theorem} %

\noindent In the case $r = 1$ the Theorem \ref{desing} reproduces an
announced result of Ciocan-Fontanine-Maulik-Kim, see
\cite{ciocan:hilb}.  Results of Konvalinka-Ciocan-Fontanine-Pak
\cite{ciocan:hilb} substantially simplify the ``mirror map''
$$\exp( \tau_{X,T,-}^{(1)}/\zeta) = (1 + q)^{k(\xi_- +\xi_+)/\zeta} .$$  
In the case $r > 1$ the factor $\exp( \tau_{X,T,-}^{(1)}/\zeta)$ is
the identity.  Formulas for quantum multiplication on these moduli
spaces, and connections with integrable systems, are developed in
Maulik-Okounkov \cite{mo:qgqc}.

\section{Gauged Gromov-Witten invariants} 

In this section we review the construction of gauged Gromov-Witten
invariants.

\subsection{Mundet stability} 
\label{stab}

Mundet stability combines the slope conditions from Ramanathan
stability for bundles and Hilbert-Mumford stability for points in the
target.  First we recall Mumford-Seshadri stability.  Let $C$ be a
smooth projective curve and $E \to C$ a vector bundle of vanishing
degree $\deg(E) = (c_1(E),[C])$.  The bundle
$$ E \ \text{{\em semistable} resp.  {\em stable}} \iff ( \deg(F) \leq
0 \quad \text{resp}. < 0, \quad \forall F \subset E) $$
for all holomorphic sub-bundles $F \subset E$ \cite{ns:st}.  In the case of a
rational curve the Birkhoff-Grothendieck theorem \cite{gr:cl} shows
that any bundle splits as a sum of line bundles.  Semistability in
degree zero on rational curves is simply the condition that the bundle
is trivial and there are no stable bundles.

Ramanathan's stability \cite{ra:th} generalizes the Mumford-Seshadri
condition to principal bundles as a condition on parabolic reductions.
Let $G$ be a connected reductive group with Lie algebra $\g$.  Let
$T \subset G$ be a maximal torus, with Lie algebra $\t$.  Denote the
integral resp. rational weights resp. coweights
$$\t_\Z = \exp^{-1}(e), \quad \t^\dual_\Z \subset \t^\dual :=
\Hom(\t,\C), \quad \t_\Q = \t_\Z \otimes_\Z \Q, \quad \t_\Q^\dual =
\t_\Z^\dual \otimes_\Z \Q $$
As in \eqref{rdec} let $ \RR = \RR_+ \cup \RR_- \subset
\t^\dual_\Z $ denote a set of positive and negative roots so that
$$ \g \cong \t \oplus \bigoplus_{\alpha \in \RR_-} \g_\alpha \oplus
\bigoplus_{\alpha \in \RR_+} \g_\alpha .$$
A {\em parabolic subgroup} of $G$ is a subgroup $Q$ such that $G/Q$ is
complete.  Up to conjugacy this means that the Lie algebra $\q$ of $Q$
is given by
$$ \q = 
\t \oplus \bigoplus_{ \alpha \in \RR_-} \g_\alpha \oplus
\bigoplus_{ \alpha \in \RR_Q} \g_\alpha $$
for some subset of the roots $\RR_Q \subset \RR_+$ such that $\q$ is a
Lie subalgebra of $\g$.  A {\em Levi subgroup} of $Q$ is a maximal
reductive subgroup $L(Q)$; again up to conjugacy the Lie algebra
$\l(\q)$ of $L(Q)$ resp. $\u(q)$ of a maximal unipotent $U(Q)$ is
$$ \l(\q) = \t \oplus \bigoplus_{ \alpha \in -\RR(Q)} \g_\alpha \oplus
\bigoplus_{ \alpha \in \RR(Q)} \g_\alpha, \quad \u(\q) = \bigoplus_{
  \alpha \in -\RR(Q)} \g_\alpha .$$
The parabolic subgroup and its Lie algebra admit decompositions into
reductive and unipotent parts
$$ \q = \l(\q) \oplus \u(\q), \quad Q = L(Q) U(Q) .$$
Taking the quotient by the maximal unipotent gives a projection 
$$ \pi_Q : Q \to Q/U(Q) \cong L(Q) .$$
This projection has the following alternative description.  A {\em
  dominant coweight} for $Q$ is a coweight $\lambda \in \t$ such that
$$ (  \alpha (\lambda) \ge 0 , \quad \forall \alpha \in \RR_+)  \quad
\text{and} \quad (\alpha(\lambda) = 0 , \quad \forall \alpha \in \RR(Q)) .$$
Any rational coweight for $Q$ determines a one-parameter subgroup
$$\phi_\lambda: \C^\times \to Q ,\quad z \mapsto \phi_\lambda(z).$$
If $\lambda \in \q$ is a dominant rational coweight then
$$ \pi_Q(q) = \lim_{z \to 0 } \Ad(\phi_{-\lambda}(z)) q .$$
Choose an equivariant identification $\g \to \g^\dual$ that identifies
the subspaces of rational weights and coweights $\t_\Q \to
\t_\Q^\dual$.  The identification $\lambda \in \t$ determines a
rational weight $\lambda^\dual \in \t^\dual$.  After finite cover
$\lambda$ defines a one-dimensional representation
$$ \chi_{\lambda^\dual}: Q \to \C^\times, q \mapsto \chi_{\lambda^\dual}(q) $$
which factors through $L(Q)$.  

The analog for principal bundles of the stability condition for
sub-bundles is a condition for parabolic reductions together with
dominant coweights.  Let $P \to C$ be a principal $G$ bundle on a
curve $C$ over a scheme $S$; bundles are by assumption locally trivial
in the \'etale topology.  A {\em parabolic reduction} is a section
$\sigma: C \to P/Q$.  Any parabolic reduction induces a reduction of
structure group given by a sub-bundle $\sigma^*(P) \subset P$ with
structure group $Q$, given by pull-back of the $Q$-bundle $P \to P/Q$.
We denote by
$$\Gr(P) := \pi_{Q,*} \sigma^*P \to C $$
the corresponding $L(Q)$-bundle, called the {\em associated graded}
bundle of $P$ associated to the parabolic reduction $\sigma$.  In case
$G = GL(n)$, a parabolic reduction is equivalent to a partial flag of
sub-bundles $E^{i_1} \subset E^{i_2} \subset \ldots \subset E^{i_l} =
E$ in the associated vector bundle $E = P(\C^n)$; the corresponding
parabolic reduction $\sigma^*P$ is the bundle of frames whose first
$i_k$-elements belong to $E^{i_k}$ for $k = 1,\ldots, l$.  The
associated graded principal bundle is the principal bundle of frames
of the associated graded vector bundle
$$ \Gr(E) = \oplus_j (E^{i_{j+1}}/E^{i_j}), \quad \Gr(P) = \Fr(\Gr(E))
.$$
\noindent The construction of the associated graded bundle also has an
interpretation via degeneration.  The family of elements
$\phi_\lambda(z)$ defines a family of automorphism
$ \Ad(\phi_\lambda(z)) : G \to G $.  Consider the family of bundles
$P^\lambda \to C \times \C^\times$ by conjugating the transition maps
of $\sigma^* P$ by $\phi_\lambda(z)^{-1}$.  Then $P^\lambda$ extends
over the central fiber $C \times \{ 0 \}$ as the bundle $\Gr(P)$.  The
{\em Ramanathan weight} of a principal bundle with respect to a
parabolic reduction and dominant weight is the degree of the line
bundle corresponding to the given dominant coweight:
$$ \mu(\sigma,\lambda) = \deg( \pi_{Q,*} \sigma^*P \times_{L(Q)}
\C_\lambda) = ([C], c_1 \pi_{Q,*} \sigma^*P \times_{L(Q)} \C_\lambda)
.$$
Then  
$$ P \ \text{ {\em semistable} resp. {\em stable}} \iff
\mu(\sigma,\lambda) \leq 0 \quad \text{resp.} < 0, \forall
(\sigma,\lambda) .$$
For rational curves Birkhoff-Grothendieck \cite{gr:cl} again implies
that a principal bundle with vanishing degree is semistable iff it is
trivial.  As for vector bundles, it suffices to check the condition
for reduction to {\em maximal} parabolic subgroups $Q$.  Ramanathan
\cite{ra:th} shows the existence of a projective coarse moduli space
for semistable principal bundles with reductive structure group and
fixed numerical invariants.

Mundet semistability \cite{mund:corr, schmitt:univ} generalizes
Ramanathan stability to the case of maps to a quotient stack.  Let $G$
be a connected reductive group acting on a smooth projective variety
$X$.  By a {\em gauged map} with domain a curve $C$ we mean a map from
$C$ to the quotient stack $X/G$, given by a pair $(P,u)$ of a
$G$-bundle and section of the associated $X$-fiber bundle:
$$ P \to C, \quad u: C \to P \times_G X .$$
Given a pair $(P \to C, u : C \to P(X))$, the section $u$ defines a
section $u^\lambda$ of $P^\lambda$ as follows: In any local
trivialization $P(X)|U \cong U \times X$ the section $u$ is given by a
map $u|U : U \to X$, and the sections $\phi_\lambda(z) u$ patch
together to a section of $P^\lambda(X)$.  By Gromov compactness,
$u^\lambda$ extends over the central fiber $C \times \{ 0 \}$ as a
stable map denoted $\Gr(u): \hat{C} \to \Gr(P)(X)$.  Associated to
this limit there is an associated {\em Hilbert-Mumford weight} defined
as follows.  The principal component $C_0$ of $\hat{C}$ is the
irreducible component such that the restriction $u_0$ of $u$ to $C_0$
maps isomorphically to $C$.  The principal component $\Gr(u)_0$ of the
associated graded section $\Gr(u)$ takes values in the fixed point set
$(\Gr(P)(X))^\lambda = \Gr(P)(X^\lambda)$ of the infinitesimal
automorphism of $\Gr(P)(X)$ induced by $\lambda$.  The {\em
  Hilbert-Mumford weight}
\begin{equation} 
\mu_H(\sigma,\lambda) \in \Z \end{equation} 
determined by the linearization $\ti{X}$, is the weight of the
$\C^\times$-action generated by $-\lambda$ on the fiber of the bundle
$(\Gr(P))(\ti{X}) \to (\Gr(P))(X)$ over a generic value of $\Gr(u)_0$:
$$ \phi_\lambda(z) \ti{x} = z^{\mu_H(\sigma,\lambda)} \ti{x}, \quad z
\in \C^\times .$$
The {\em Mundet weight} is the sum of the Hilbert-Mumford and
Ramanathan weights:
$$ \mu_M(\sigma,\lambda) : = \mu_H(\sigma,\lambda) +
\mu_R(\sigma,\lambda). $$
Then 
$$ (P,u) \ \text{ {\em semistable} resp. {\em stable}} \iff
\mu(\sigma,\lambda) \leq 0 \quad \text{resp.} < 0, \forall
(\sigma,\lambda) .$$
Mundet's original definition allowed possibly irrational $\lambda$,
but this is unnecessary in the case that the symplectic class is
rational by \cite[Remark 5.8]{qk2}.  Mundet semistability is realized
as a git stability condition in Schmitt \cite{schmitt:univ,
  schmitt:git}.

The moduli stack of Mundet-semistable morphisms admits a natural
Kontsevich-style compactification that allows formation of bubbles in
the fibers of the associated bundle: An {\em $n$-marked gauged map}
from ${C}$ to $X$ over a scheme $S$ is a datum $(\hat{C},P,u,\ul{z})$
where $\hat{C} \to S$ is a proper flat morphism with reduced nodal
curves as fibers, $ P \to {C} \times S$ is a principal $G$-bundle; and
$${u}: \hat{C} \to P(X) := (P \times X)/G$$ 
is a family of stable maps with base class $[{C}]$, that is, the
composition of ${u}$ with the projection $P(X) \to {C}$ has class
$[{C}]$.  A {\em morphism} between gauged maps $(S,\hat{C},P,u)$ and
$(S',\hat{C}',P',u')$ consists of a morphism $\beta: S \to S'$, a
morphism $\phi: P \to (\beta \times 1)^*P'$, and a morphism $\psi:
\hat{C} \to \hat{C}'$ such that the first diagram below is Cartesian
and the second and third commute:
$$ \begin{diagram} \node{\hat{C}} \arrow{e}
  \arrow{s,l}{\psi} \node{S} \arrow{s,r}{\beta} \\
  \node{\hat{C}'} \arrow{e} \node{S'} \end{diagram} \quad
\begin{diagram} \node{P} \arrow{e} \arrow{s,t}{\phi} \node{S
	\times {C}} \arrow{s,b}{\on{id}} \\ \node{(\beta \times
	1)^*P'} \arrow{e} \node{S \times {C}} \end{diagram}
  \quad \begin{diagram} \node{\hat{C}} \arrow{e,t}{u}
	\arrow{s,t}{\psi} \node{P(X)} \arrow{s,b}{[\phi \times
	\on{id_X}]} \\ \node{\hat{C}'} \arrow{e,b}{u'}
	\node{P'(X).} \end{diagram} $$
An {\em $n$-marked} nodal gauged map is equipped with an $n$-tuple
$(z_1,\ldots,z_n) \in \hat{C}^n$ of distinct smooth points on
$\hat{C}$.  An $n$-marked nodal gauged map $(\hat{C},P,\ul{z},u)$ is
{\em Mundet semistable} resp. {\em stable} if the principal component
is Mundet semistable resp. stable and the section $u: \hat{C} \to
P(X)$ is a stable section, in the sense that any component on which
$u$ is constant has at least three special (nodal or marked) points.

\subsection{Moduli stacks} 

We introduce the following notations for moduli stacks.  Denote by
$\ol{\MM}_n^G(C,X,d)$ resp. $\ol{\M}_n^G(C,X,d)$ the category of
gauged maps resp. Mundet semistable gauged maps from $C$ to $X/G$ of
homology class $d$ and $n$ markings.

\begin{theorem} 
For any $d,n$, if stable=semistable then the stack
$\ol{\M}_n^G(C,X,d)$ is a proper Deligne-Mumford stack equipped with
evaluation morphisms
$$\ev: \ol{\M}_n^G(C,X,d) \to (X/G)^n, \quad (\hat{C},P,u) \mapsto
(\ul{z}^*P, \ul{z}^* u) $$
and virtual fundamental class. 
\end{theorem} 

The properties of the moduli stacks in the above theorem were proved
elsewhere.  Properness is covered in detailed in \cite[Theorem
1.1]{reduc}.  Virtual fundamental classes are \cite[Example 6.6]{qk2}.
We sketch the construction for completeness.  The proof of properness
uses a simpler Grothendieck-style compactification obtained by
allowing the maps to acquire base points, studied by Schmitt
\cite{schmitt:univ}, \cite[Section 2.7]{schmitt:git}.  Suppose that
$X \subset \P(V)$ is embedded in the projectivization $\P(V)$ of a
$G$-representation $V$.  A map $C \to P(\P(V))$ gives rise to a line
sub-bundle $L \subset C \times P(V)$.  By dualization such a
sub-bundle gives rise to a quotient map
$q: C \times P(V)^\dual \to L^\dual$.  A {\em gauged quotient} is a
datum $(P,L,q,\ul{z})$, called by Schmitt \cite{schmitt:univ} a {\em
  bundle with map}. Denote by $\ol{\M}_n^{G,\quot}(C,X,d)$ the space
of stable gauged quotients.  The moduli stacks
$\ol{\M}_n^{G,\quot}(C,X,d)$ only admit evaluation morphisms to the
quotient stacks for the ambient vector spaces,
$$ \ev: \ol{\M}_n^G(C,X,d) \to (V/(G \times \C^\times))^n, \quad
(\hat{C},P,u) \mapsto (\ul{z}^*P, \ul{z}^* L, \ul{z}^*q) .$$

The moduli stack of stable gauged quotients admits a construction as a
geometric invariant theory quotient by Schmitt
\cite{schmitt:univ,schmitt:git}.  Choose a faithful representation
$G \to GL(V)$, so that $X \subset \P(V)$.  A {\em $k$-level structure}
for a stable gauged quotient is a collection of sections
$s_1,\ldots,s_k: C \to P(V) $ generating $P(V)$.  Equivalently, a
level structure is a surjective morphism
$ \mO_C^{\otimes k} \to P(V)^\dual $.  The action of $GL(k)$ on
$C \times \C^k$ induces an action on the set of level structures by
composition.  The stack $\ol{\M}_n^{G,\lev,\quot}(C,X,d)$ of gauged
quotients with level structure is naturally an Artin stack with an
action of $GL(k)$ on the sections.  Schmitt \cite[Section
2.7]{schmitt:univ} constructs a linearization
$D(\ti{X}) \to \ol{\M}_n^{G,\lev,\quot}(C,X,d)$ giving rise to a
projective embedding of the coarse moduli space, so that the git
quotient is the stack of gauged quotients:
$$ \ol{\M}_n^{G,\quot}(C,X,d) = \ol{\M}_n^{G,\lev,\quot}(C,X,d) \qu
GL(k) .$$
In particular this construction implies that
$\ol{\M}_n^{G,\quot}(C,X,d)$ has proper coarse moduli space.  If
stable=semistable then all stabilizers are finite, and since we are in
characteristic zero, this implies that $\ol{\M}_n^{G,\quot}(C,X,d)$ is
Deligne-Mumford and proper.  Now the Kontsevich-style compactification
$\ol{\M}_n^G(C,X,d)$ admits a morphism by Popa-Roth \cite[Theorem
  7.1]{po:stable}
$$ \ol{\M}_n^G(C,X,d) \to \ol{\M}_n^{G,\quot}(C,X,d) $$
and so is also proper.  Denote by $ \ol{\M}_n^{G,\lev}(C,X,d)$ the
moduli stack of gauged maps with level structure on the associated
vector bundle $P(V)$.  The Givental construction on the moduli stack
of maps with level structure gives a morphism $\pi:
\ol{\M}_n^{G,\lev}(C,X,d) \to \ol{\M}_n^{G,\lev,\quot}(C,X,d)$.  Then
the moduli stack of gauged maps is also a stack-theoretic quotient
$$ \ol{\M}_n^{G}(C,X,d) = \pi^{-1}(\ol{\M}_n^{G,\lev}(C,X,d)^{\ss}) / G .$$
However, the pull-back of the linearization $D(\ti{X})$ is not ample
on $\ol{\M}_n^{G,\lev}(C,X,d)$.  Thus this quotient cannot be
considered a git quotient without further perturbation of the
linearization.

Virtual fundamental classes are obtained from the construction of
Behrend-Fantechi \cite{bf:in}.  The argument uses the deformation
theory from Olsson \cite[Theorem 1.5]{ol:def} for morphisms to
quotient stacks.  The universal curve $\ol{\cC}_n^G(C,X)$ is the stack
whose objects are tuples $(\hat{C},P,u,\ul{z},z')$ where
$(\hat{C},P,u,\ul{z})$ is a gauged map and $z' \in \hat{C}$ is a
(possibly singular) point.  Forgetting $z'$ defines a projection
$$ p:\ol{\cC}_n^G(C,X) \to \ol{\M}_n^G(C,X) $$
while evaluating at $z'$ defines a universal gauged map
$$ e: \ol{\cC}_n^G(C,X) \to X/G .$$
The relative obstruction theory has complex given by $ R\pi_* e^*
T(X/G)^\dual$ equipped with its canonical morphism to the cotangent
complex of $\ol{\M}_n^G(C,X)$. If stable=semistable then the
obstruction theory is perfect and $\ol{\M}_n^G(C,X,d)$ is a proper
smooth Deligne-Mumford stack with perfect relative obstruction theory
over the stack of {\em semistable} $n$-marked maps to $C$, see
\cite{qk2}.  Denote by $ [\ol{\M}_n^G(C,X,d)] \in H(\ol{\M}_n^G(C,X,d)) $ the
virtual fundamental classes constructed via Behrend-Fantechi
machinery.  

Using the virtual fundamental classes, gauged Gromov-Witten potentials
are defined as follows.  Suppose that stable=semistable for all gauged
maps. The gauged potential $\tau_{X,G}$ is the formal map defined by
$$ \tau_{X,G}: QH_G(X) \to \Lambda_X^G $$
$$ 
{h} \mapsto \sum_{n \ge 0, d \in H_2^G(X,\Z)/\on{torsion}}
(q^d/n! ) \int_{ [\ol{\M}_n^G(C,X,d)] } \ev^* ({h},\ldots, {h}
) $$
for ${h} \in H_G(X)$.

Later we will need several variations on the gauged Gromov-Witten
potential.  We describe three variations which will be used later.

\begin{definition} 
\begin{enumerate} 
\item {\rm (Gauged invariants with Deligne-Mumford classes)} The first
  variation involves pull-back classes from the curve $C$.  That is,
  let
$$ f: \ol{\M}_n^G(C,X,d) \to \ol{\M}_n(C) $$
be the map obtained by projecting from $C \times X/G$ to $C$.  For any
class $\beta \in H(\ol{\M}_n(C))$ and ${h} \in H_G(X)^n$ define a
gauged invariants with insertions
\begin{equation} \label{ins}
 \tau_{X,G}^{n,d}({h},\beta) := \int_{ [\ol{\M}_n^G(C,X,d)] } \ev^* (h
 \otimes \ldots \otimes h) \cup f^* \beta.\end{equation}
In particular, by taking $n = 3$ and $\beta$ the dual class of a point
in $\ol{\M}_3(C)$ we obtain the gauged analog of three-point
invariants.
\item {\rm (Twisted invariants)} The second variation gives
  Euler-twisted gauged Gromov-Witten invariants. For any
  $G$-equivariant bundle $E \to X$ we denote by
\begin{equation} \label{index} \Ind(E) :=  Rp_* e^* (E/G)
\end{equation} 
the index of the bundle $E/G \to X/G$.  The index class $\Ind(E)$ lies
in the bounded derived category of $\ol{\M}_n^G(C,X)$, since $p$ is
proper. 
 Furthermore $\Ind(E)$ admits a resolution by vector bundles,
since $p$ is a local complete intersection morphism, see
\cite[Appendix]{co:qrr}.  It follows that the Euler class
\begin{equation} \label{epsE} \eps(E) :=
  \Eul_{\C^\times}(\Ind(E)) \in
  H_{\C^\times}(\ol{\M}_n^G(C,X)) \end{equation}
is well-defined after passing to the equivariant cohomology of
$\ol{\M}_n^G(C,X)$ for the trivial $\C^\times$-action corresponding to
scalar multiplication on the fibers and inverting the equivariant
parameter.  The Euler-twisted gauged invariants are defined by
\begin{multline} \label{twisted} \tau_{X,G}^{n,d}: QH_G(X)^n
  \times H(\ol{\M}_n(C)) \to \Q \\ ({h},\beta) \mapsto \int_{
    [\ol{\M}_n^G(C,X,d)]} \ev^* ({h},\ldots, {h} ) \cup f^* \beta \cup
  \eps(E) .\end{multline}
\item {\rm (Parabolic structures)} A final variation involves adding a
  parabolic structure at a point on the curve as in Heinloth-Schmitt
  \cite{hs:coh} and Beck \cite{beck:prin}.  This variation will be
  used later to shift the stability condition slightly so that certain
  stacks become Deligne-Mumford.  Recall from e.g. \cite{bs:mo} that a
  {\em quasi-Borel structure} on a $G$-bundle $P \to C$ at $z_0 \in C$
  consists of a reduction of structure group $ \rho_{z_0} \subset
  P_{z_0}/B $.  In the case of $GL_r$, a quasi-Borel structure is a
  full flag in the associated vector bundle $E = P(\C^r)$ at the point
  $z_0.$ A {\em Borel-structure} is a quasi-Borel structure
  $\rho_{z_0}$ together with a {\em dominant weight} $\nu \in
  \t_\Q^\dual$ for $B$ that lies in the interior of the Weyl alcove;
  That is, we have $(\alpha_0,\nu) < 1$ with respect to the basic
  inner product on $\t^\dual$.  We consider here only generic small
  parabolic weights, see \cite{bs:mo} for the general theory.

The definition of the Ramanathan weight extends to a Ramanathan weight
for bundles with parabolic structure, with an additional term arising
from the parabolic structure \cite{bs:mo}.  The parabolic weight
defines a line bundle over the generalized flag variety $ G \times_B
\C_\nu \to G/B .$ Given any point $gB \in G/B$, the limit of $gB$
under $z^\lambda, z \to 0$ is determined by the Bruhat cell $ B w B
\ni g B $ containing $g$.  We denote by $\mu_B(\lambda)$ the
corresponding Hilbert-Mumford weight,
$$ \mu_B(\lambda) = (w \nu, \lambda) \in \Q.$$
Now let $P^\lambda \to C \times \C^\times$ denote the family of
bundles in Section \ref{stab}.  The reduction $\rho_{z_0}$
defines a reduction $\rho_{z_0}$ in $P^\lambda$ and, since $G/B$ is
complete, extends over the central fiber to a reduction
$\Gr(\rho)_{z_0} \in \Gr(P)_{z_0}/B$.  The automorphism $\lambda$ acts
on $\Gr(P_{z_0}) \times_B \C_\mu$ by some integer
$\mu_B(\sigma,\lambda)$.  The Ramanathan weight becomes modified by
the addition
$$ \mu_R(\sigma,\lambda,\nu) = \mu_R(\sigma,\lambda) +
\mu_B(\sigma,\lambda) = \mu_R(\sigma,\lambda) + 
(w\nu,\lambda) $$
where $w$ is the Weyl group element corresponding to the Bruhat cell
containing $\sigma(z_0)$.  

We introduce the following notations for moduli stacks of gauged maps
with parabolic structure.  Let $\ol{\M}_n^G(C,X,\ti{X},\nu)$
denote the Artin stack consisting of Mundet semistable pairs
$(P,u,\sigma)$ of a bundle $P \to C$, a section $u: C \to P(X)$,
markings $\ul{z}$, and a parabolic structure $(\sigma,\nu)$ at $z_0$.
\end{enumerate} 
\end{definition} 

\begin{proposition} \label{irrat} 
 If stable=semistable then $\ol{\M}_n^G(C,X,\ti{X},\nu)$ is a smooth,
 proper Deligne-Mumford stack with a perfect relative obstruction
 theory.  For generic parabolic weights $\nu$, stable=semistable and
 so the conclusion holds.  If stable=semistable for
 $\ol{\M}_n^G(C,X,\ti{X})$ and $\nu$ is sufficiently small, then the
 parabolic structure does not play a role in stability and the
 forgetful morphism
$$\pi: \ol{\M}_n^G(C,X,\ti{X},\nu) \to \ol{\M}_n^G(C,X,\ti{X})$$
is a $G/B$-bundle.
\end{proposition} 

\begin{proof}[Sketch of proof]  
  By considering a principal bundle as a vector bundle together with
  section of the associated $GL(r)/G$ bundle as in \cite{beck:prin},
  the construction of the quot-scheme compactification of the moduli
  stacks reduces to the construction of vector bundles with local and
  global decorations in \cite{beck:dec}.  On the other hand, the
  Kontsevich-style compactification $\ol{\M}_n^G(C,X,\ti{X},\nu)$ is
  obtained from the quot-scheme compactification by taking stable
  sections of the associated fiber bundle.  It follows that
  $\ol{\M}_n^G(C,X,\ti{X},\nu)$ is an Artin stack, and Deligne-Mumford
  if stable=semistable.  For generic weights $\nu$ equality cannot
  hold in the semistability inequality.  It follows that
  stable=semistable.  The construction of the moduli space of gauged
  maps with parabolic structure is a standard extension of the
  construction of moduli of parabolic bundles, and omitted.  The final
  assertion follows from the fact that for small $\nu$, the parabolic
  structure plays no role in the stability condition. Forgetting the
  parabolic structure gives a morphism to the moduli stack
  $\ol{\M}_n^G(C,X,\ti{X})$ with fiber $P_{z_0}/B \cong G/B$.
\end{proof} 

\begin{remark}\label{upstairs} 
 The following construction will be used later to avoid singularities
 in the ``master space'' used for wall-crossing.  If stable=semistable
 then the integral of any class ${h}$ over $\ol{\M}_n^G(C,X,\ti{X})$
 may be written as an integral over the moduli space of maps with
 parabolic structure: Let $L_\pi$ be the relative cotangent complex
 for the projection $\pi: \ol{\M}_n^G(C,X,\ti{X},\nu) \to
 \ol{\M}_n^G(C,X,\ti{X})$.  Then
\begin{equation} \label{fiberint}  \int_{[\ol{\M}_n^G(C,X,\ti{X})]} {h} = 
\int_{[\ol{\M}_n^G(C,X,\ti{X},\nu)]}
 \pi^* {h} \cup \Eul(L_\pi)/ |W| .\end{equation} 
Indeed the integral of $\Eul(L_\pi)/|W|$ over the fiber of $\pi$ is
$$ \int_{[G/B]} \Eul(T^\dual(G/B))/|W| = \chi(G/B)/|W| = 1 .$$
\end{remark} 

\subsection{Localized graph potentials} 

Restricting to the case of a rational curve, one may factorize the
graph potential into {\em localized gauged potentials} corresponding
to the two fixed points:

\begin{theorem}  \cite{gi:eq}, \cite{qk3}  There exist {\em localized graph resp. gauged potentials}
$$ \tau_{X \qu G,\pm}: QH(X \qu G) \to QH_{\C^\times}(X \qu G), \quad
  \tau_{X,G,\pm}: QH_G(X) \to QH_{\C^\times}(X \qu G),
 $$
that represent the contribution from the fixed points $0$
resp. $\infty$, such that the graph resp. gauged potential then admits
a factorization
\begin{multline} \tau_{X \qu G} = (\tau_{X \qu G,\pm}, \tau_{X \qu G,+}) : QH(X \qu
G) \to \Lambda_X^G \\ \tau_{X,G} = (\tau_{X,G,\pm}, \tau_{X,G,+}) :
QH_G(X) \to \Lambda_X^G \end{multline}
where $( \cdot , \cdot )$ denotes the pairing given by integration
over $I_{X \qu G}$.
\end{theorem} 

The proof is given elsewhere \cite{qk3}.  We sketch the construction
for completeness.  Let $\C^\times$ act on $C = \P^1$ via the standard
action with weights $-1,1$:
$$ \C^\times \times \P^1 \to \P^1 , \quad (w,[z_-,z_+]) \mapsto
w[z_-,z_+] = [w^{-1}z_-, w z_+] .$$
The fixed points of the induced action on $\M_n(C,X \qu G)$ correspond
to configurations of a constant map to $X \qu G$ together with bubble
trees attached at the fixed points $[1,0],[0,1]$.  That is, the fixed
point locus is a union of fibered products
$$ \M_n(C,X \qu G)^{\C^\times} = \bigcup_{n_- + n_+ = n} \bigcup_{d_ +
  d_+ = d} \ol{\M}_{0,n_-+1}(X \qu G,d_-) \times_{\ol{I}_{X \qu G}}
\ol{\M}_{0,n_++1}(X \qu G,d_+) .$$
By pushing-forward the classes $\ev^* (h \otimes \ldots \otimes h)$ over the
extra marked point in $\ol{\M}_{0,n_\pm + 1}(X \qu G,d_\pm)$ define
$$ \tau_{X \qu G,\pm}: QH(X \qu G) \to QH_{\C^\times}(X \qu G) $$
that represent the contribution to $\tau_{X \qu G}$ from the fixed
points $0$ resp. $\infty$.   More precisely,
$$ \tau_{X \qu G,\pm}({h}) := 1 + h/\zeta + \sum_{n \ge 0} (1/n!)
\tau_{X \qu G,\pm}^n(h \otimes \ldots \otimes h) $$
where 
$$ \tau_{X \qu G,\pm}^n({h}_1,\ldots,{h}_n) = \sum_{d \in H_2(X
  \qu G,\Z)} q^d/n! \ev_{n+1,*} \left( \mp \zeta (\pm \zeta -
\psi_{n+1})^{-1} \bigcup_{i=1}^n \ev_i^* {h}_i \right)$$
is the sum over $n,d$ such that the moduli space of stable maps is
non-empty (that is, either $ n \ge 3$ or $d > 0$) and the class
$$\psi_{n+1} = c_1(L_{n+1}) \in H^2(\ol{\M}_{0,n+1}(X \qu G,d))$$ 
is the first Chern class of the cotangent line at the $(n+1)$-st
marked point,
$$ L_{n+1} \to \ol{\M}_{0,n+1}(X \qu G,d)), \quad (L_{n+1})_{ (u: C
  \to X \qu G, \ul{z})} = T_{z_{n+1}} C .$$
The graph potential then admits a factorization
$$ \tau_{X \qu G} = (\tau_{X \qu G,\pm}, \tau_{X \qu G,+}) : QH(X \qu G) \to \Lambda_X^G $$
where $( \cdot , \cdot )$ denotes the pairing given by integration
over $I_{X \qu G}$ \cite{gi:eq}.  We denote, for later use, the map
\begin{equation} \label{infty} 
 \tau_{X \qu G,\pm}^\infty: QH(X \qu G)^2 \to \Lambda_X^G, \quad
 \tau_{X \qu G,\pm}^\infty(h_1,h_2) = (\tau_{X \qu G,\pm}(h_1),
 h_2) \end{equation} 
obtained by dualizing one factor.

The fixed points for the circle action on the space of gauged maps
over the projective line are described in \cite{qk3}: Fixed maps are
data $(\hat{C},P,u,\ul{z})$ such that there exists a one-parameter
family of bundle isomorphisms preserving the section:
$$ \phi: \C^\times \to \Hom(P, m_w^*P) \quad P(\phi(w))^* m_w^*u = u,
\forall w \in \C^\times $$
where $m_w: C \to C$ is multiplication by $C$.  In local
trivializations near the fixed points $[1,0],[0,1]$, the bundle
automorphism $\phi$ is given by homomorphisms
$$ \phi_\pm : \C^\times \to G .$$
In the corresponding trivializations of $P(X)$, the section $u$ is given by maps
$$ u_\pm: \C \to G, \quad \phi_\pm(w) u_\pm(z) = u_\pm(wz) .$$
Furthermore, for $u: \hat{C} \to P(X)$ a $\C^\times$-fixed map sends
the components of $\hat{C}$ map to either of the fixed points $[1,0],
[0,1] \in (\P^1)^{\C^\times} $.  Thus $u$ consists of a pair of
sections $(u_-,u_+)$ with bubble trees attached at $[1,0],[0,1] \in
\P^1$, glued together via a transition map 
$$\C^\times \to G, \quad z \mapsto \phi_+(z) \phi_-(z^{-1})^{-1} .$$

One may therefore view the fixed point locus as a fiber product, as
follows: A {\em framing} of a gauged map is a trivialization of $P$
neighborhood of a point $z' \in C$.  We consider the stack
$\ol{\M}^{G,\fr}_{n_\pm,\pm}(C,X,d)^{\C^\times}$ of $\C^\times$-fixed
{\em framed gauged maps} at the point $z' =[1,0]$ resp. $z' = [0,1]$.
In the case of large linearization, evaluation at the extra marked
point defines maps
$$ \ev': \ol{\M}^{G,\fr}_{n_\pm,\pm}(C,X,d) \to I_{X \qu G} .$$
Then the fixed point locus factorizes 
$$ \ol{\M}_n^{G}(C,X,d)^{\C^\times} =
\ol{\M}^{G,\fr}_{n_-,-}(C,X,d_-)^{\C^\times} \times_{I_{X \qu G}}
\ol{\M}^{G,\fr}_{n_+,+}(C,X,d_+)^{\C^\times} .$$
By pushing-forward classes $\ev^* (h \otimes \ldots \otimes h)$ we obtain maps
$$ \tau_{X,G,\pm}: QH_{G}(X) \to QH_{\C^\times}(X \qu G) $$
so that 
$$ \tau_{X,G} = (\tau_{X,G,+}, \tau_{X,G,-}) : QH_G(X) \to \Lambda_X^G .$$

Later we will need a variation on this construction that involves
dualizing one of the factors.  For any real parameter $\rho$ let
$\ol{\M}^G_n(C,X,\ti{X}^\rho,d)_- \subset
\ol{\M}^G_n(C,X,\ti{X}^\rho,d)$
denote the locus of the fixed point set where all markings
$z_1,\ldots, z_n$ map to $0 \in C$; all bubble components
$C_1,\ldots, C_k$ map to $0 \in C$; the one-parameter subgroup
$\phi_+$ vanishes.  Evaluation at $z_1,\ldots, z_n$ and at
$\infty \in C$ defines a map
$$ \ev_- \times \ev_+: \ol{\M}^G_n(C,X,\ti{X}^\rho,d)_- 
\to (X/G)^n
\times X/G .$$
By integration one obtains a map linear in the second variable
\begin{multline} \tau^\rho_{X,G,-}:  QH_G(X)^2 \to \Lambda_X^G[\zeta,\zeta^{-1}]],  \\
(h_-, h_+) \int_{[\ol{\M}^G_n(C,X,\ti{X}^\rho,d)_-]} \left(
    \ev_-^*(h_- \otimes \ldots \otimes h_- \right) \cup \ev^*_+ h_+
    \cup \Eul(\nu_-) \end{multline}
where $\nu_-$ is the virtual normal complex of
$\ol{\M}^G_n(C,X,\ti{X}^\rho,d)_- $ in
$\ol{\M}^G_n(C,X,\ti{X}^\rho,d)$.  In the limit $\rho \to \infty$, the
gauged map is generically semistable and triviality in a neighborhood
of $\infty$ implies that in fact $\ev_+$ takes values in $X \qu G$.
Hence in this case $\tau^\rho_{X,G,-}$ is the obtained from the
dualization
$$\tau_{X,G,-}^\dual : QH_G(X) \times QH(X \qu G) \to \Lambda_X^G$$
of $\tau_{X,G,-}$ by composition with pull-back $QH_G(X) \to QH(X \qu
G)$.  Similar definitions requiring the bubbles and markings to map to
$\infty$ give rise to a map $ \tau^\rho_{X,G,+}: QH_G(X)^2 \to
\Lambda_X^G$ related to $\tau_{X,G,+}$ in the limit $\rho \to \infty$
by dualization.

\subsection{Master space}

In this section we study the area-dependence of the gauged
Gromov-Witten invariants, by which we mean the dependence on the
choice of linearization.  The basic strategy is the same as that
outline in e.g. Thaddeus \cite{th:fl} for the case of variation of
linearization in geometric invariant theory.  The wall-crossing
formula of is obtained from a {\em master space} construction as
follows: Suppose that
$$\ti{X}_\pm \to X$$
are linearizations, that is, ample $G$-line bundles on $X$.
Associated to $\ti{X}_\pm$ are linearizations
$$D(\ti{X}_\pm) \to \ol{\M}^{G,\lev,\quot}_n(C,X,d)$$ 
over the moduli stacks of curves with level structure constructed in
Schmitt \cite{schmitt:univ,schmitt:git}.  Consider the rank two bundle
obtained from the direct sum:
$$ D(\ti{X}_-) \oplus D(\ti{X}_+) \to \ol{\M}^{G,\lev,\quot}_n(C,X,d)
.$$
Taking the projectivization of the total space gives a
$\P^1$-fibration
$$ \P( D(\ti{X}_-) \oplus \ D(\ti{X}_+)) \to
\ol{\M}^{G,\lev,\quot}_n(C,X,d) .$$
The action of $GL(k)$ lifts to the fibration, since the bundles
$D(\ti{X}_\pm)$ are $GL(k)$-equivariant.  The bundle
$$\mO_{\P( D(\ti{X}_-) \oplus \ D(\ti{X}_+))}(1) \to 
\P( D(\ti{X}_-)
\oplus \ D(\ti{X}_+))$$
is automatically ample on the coarse moduli space of $\P( D(\ti{X}_-)
\oplus \ D(\ti{X}_+))$.  Denote the git quotient with respect to this
linearization
$$ \ol{\M}^{G,\lev,\quot}_n(C,X,\ti{X}_-,\ti{X}_+,d) = \P( D(\ti{X}_-)
\oplus \ D(\ti{X}_+)) \qu GL(k) .$$
Similarly, let  
$$\pi^* \P( D(\ti{X}_-) \oplus \ D(\ti{X}_+)) \to
\ol{\M}^{G,\lev}_n(C,X,d) $$
denote the pull-back to the stack of stable gauged maps with level structure and 
$$ \ol{\M}^{G,\lev}_n(C,X,\ti{X}_-,\ti{X}_+,d) = \pi^{-1}( \P(
D(\ti{X}_-) \oplus \ D(\ti{X}_+)))^{\ss} / GL(k) $$
the quotient of the pull-back of the semistable locus.  The action of
$\C^\times$ on $\P( D(\ti{X}_-) \oplus \ D(\ti{X}_+))$ induces an
action of $\C^\times$ on $\ol{\M}^{G,\lev}_n(C,X,\ti{X}_-,\ti{X}_+,d)$.

The fixed point components for the natural circle action are of two
types.  First, there there are inclusions
$$ \P( D(\ti{X}_\pm) \oplus 0 ) \to \P( D(\ti{X}_-) \oplus
D(\ti{X}_+)) $$
and isomorphisms 
$$ \P( D(\ti{X}_\pm) \oplus 0 ) \cong \ol{\M}_n^G(C,X,\ti{X}_\pm).$$
These induce embeddings
$$ \ol{\M}_n^G(C,X,\ti{X}_\pm) \to
\ol{\M}_n^G(C,X,\ti{X}_-,\ti{X}_+)^{\C^{\times}} $$
in the locus of fixed points of the $\C^\times$-action.  On the other
hand, there are fixed point components correspond to {\em reducible}
gauged maps for some stability condition interpolating between those
defined by $\ti{X}_\pm$.  Reducibility means that the fixed point
components consist of maps $v = (P,u): \hat{C} \to X/G$ that admit a
one-parameter family of automorphisms $\phi: \C^\times \to \Aut(P)$;
via evaluation at a point $\Aut(P) \to \Aut(P_z)$, any such
one-parameter family may be identified with a one-parameter family of
automorphisms of $G$ generated by some element $\lambda \in \g$.
Euler-twisted integration over the fixed point components gives rise
to fixed point contributions
$$ \tau_{X,G,\zeta,t}: QH_G(X) \to \Lambda_X^G .$$
The fixed point contributions are bubble trees consisting of maps to
the quotient stack with one-parameter automorphisms and stable maps
fixed up to isomorphism by one-parameter subgroups.  Suppose that a
gauged map $(P \to C, u: \hat{C} \to P(X))$ is reducible, that is, has
a one-parameter family of automorphism $\phi: \C^\times \to \Aut(P)$
covering the identity on the principal component so that the
associated automorphism
$$\phi(X): P(X) \to P(X), \quad \phi(X)^*
u = u .$$  
Evaluation at any fiber defines a homomorphism $\phi_z: \C^\times \to
\Aut(P_z) \cong G$ and so identifies $\phi_z$ with a one-parameter
subgroup of $G$.  Let $\lambda \in \g$ be a generator of $\phi_z$ and
$G_\lambda \subset G$ the centralizer.  The structure group of $P$
reduces to the centralizer $G_\lambda$ of $\lambda$.  Furthermore the
restriction $u | \hat{C}_0$ of $u$ to the principal component $\hat{C
}_0$ takes values in $P(X^\lambda)$ where $X^\lambda = \{x \in X |
\lambda_X(x) = 0 \}$.  Any bubble tree attached at $z \in \hat{C}_0$
must be fixed, up to isomorphism, by the action of $\phi(z) \in
\Aut(P_z(X))$.  That is, there exists a one-parameter family of
automorphisms $\psi: \C^\times \to \Aut(\hat{C})$ so that $\psi^* u =
\phi(X) \circ u$, where $\phi(X): P(X) \to P(X)$ is the automorphism
of the associated fiber bundle induced by $\phi$.

We introduce notation for these fixed point stacks and their normal
complexes as follows.
\label{fixedef}
 For each $\lambda \in \g$, let $\ol{\M}_n^G(C,X,\ti{X},\lambda)$
 denote the stack of Mundet-semistable morphisms from $C$ to
 $X/G_\lambda$ that are $\C^\times_\lambda$-fixed and take values in
 $X^\lambda$ on the principal component.  Via the inclusion $G_\lambda
 \to G$ the universal curve over $\ol{\M}_n^G(C,X,\ti{X},\lambda)$
 admits a morphism to $X/G$.  Denote by $\nu_\lambda$ the virtual
 normal complex for the morphism $\ol{\M}_n^G(C,X,\ti{X},\lambda) \to
 \ol{\M}_n^G(C,X,\ti{X})$.

The virtual fundamental classes on these fixed point stacks lead to
fixed point contributions appearing in the wall-crossing formula.  Let
$QH_{G,\fin}(X)$ denote the tensor product of $H_G(X)$ with the
sub-ring
$$\Lambda_X^{G,\fin} = \left\{ \sum_{i=1}^n c_i q^{d_i}, d_i \in
H_2^G(X), c_i \in \Q \right\} \subset \Lambda_X^G$$
of finite sums.  Let $\xi$ denote the equivariant parameter for the
action of the one-parameter subgroup generated by $\lambda$ and
$\Resid_\xi: \C[\xi] \to \C$ the residue of $\xi$ at $0$, that is, the
map taking the coefficient of $\xi^{-1}$.
\label{contrib} Virtual integration over
$\ol{\M}_n^{G}(C,X,\ti{X},\lambda)$ defines a ``fixed point
contribution''
  \begin{multline}\label{eq:contrib} \tau_{X,G,\lambda,\ti{X}}: QH_{G,\fin}(X) \to
    \ti{\Lambda}_X^G \otimes H(B \C^\times), \\ \quad {h} \mapsto
    \sum_{d \in H_{G}^2(X,\Z)} \sum_{n \ge 0} \Resid_\xi
    \int_{[\ol{\M}_n^G(C,X,\ti{X},\lambda,d)]} (q^d/n!)  \ev^* ({h},
    \ldots, {h}) \cup \Eul(\nu_\lambda)^{-1}
\end{multline}
for ${h} \in H_G(X)$.  Here we omit the restriction map $H_{G,\fin}(X)
\to H_{G_\lambda}(X)$ to simplify notation.  The following is
\cite[Theorem 3.14]{wall}.

\begin{theorem}[Wall-crossing for gauged Gromov-Witten
  potentials] \label{gwall} Let $X$ be a smooth projective
  $G$-variety.  Suppose that $\ti{X}_\pm \to X$ are linearizations
  such that semistable=stable for the stack of polarized gauged maps
  in \cite{wall}.  Then the gauged Gromov-Witten potentials are
  related by
\begin{equation} \tau_{X,\ti{X}_+,G} - \tau_{X,\ti{X}_-,G} =
  \sum_{[\lambda],t \in (-1,1)}
\frac{|W_{\lambda}|}{|W_{\C \lambda}|} 
  \tau_{X,G,\lambda,t}
\end{equation}
where the sum is over equivalence classes $[\lambda]$ of
unparametrized one-parameter subgroups generated by $\lambda
\in \g$.  
 \label{gwall2} 
Similarly the localized gauged Gromov-Witten potentials are related by
\begin{equation} \tau_{X,\ti{X}^+,G,\pm} - \tau_{X,\ti{X}^-,G,\pm} =
\sum_{[\lambda],t \in (-1,1)}
\frac{|W_{\lambda}|}{|W_{\C \lambda}|} 
 \tau_{X,\ti{X}^t,G,\pm,\lambda}: QH_G^{\fin}(X)^2 \to
\Lambda_X^G.
\end{equation}
\end{theorem} 

The fixed point contributions can be re-written as contributions from
gauged Gromov-Witten invariants with structure group of smaller rank
as follows.  For $\lambda \in \g$ let $\C^\times_\lambda \subset
G_\lambda$ denote the one-parameter subgroup generated by $\lambda$,
and $G_\lambda/\C^\times_\lambda$ the quotient.  Let $X^\lambda
\subset X$ denote the fixed point set of $\C^\times_\lambda$.  Let
$ \ol{\M}_{0,n}(X)^{\C^\times_\lambda} \subset \ol{\M}_{0,n}(X) $
denote the $\C^\times_\lambda$-fixed point stack of stable maps to
$X$.  The evaluation map restricted to
$\ol{\M}_{0,n}(X)^{\C^\times_\lambda}$ automatically takes values in
the fixed point locus $X^\lambda \subset X$, that is, $ \ev:
\ol{\M}_{0,n}(X)^{\C^\times_\lambda} \to (X^\lambda)^n .$ Push-pull
over the moduli stack $\ol{\M}_{0,n+1}(X)^{\C^\times_\lambda}$ defines
a {\em quantum restriction map}
$$ \iota_\lambda: QH_{G_\lambda}(X) \to QH_{G_\lambda}(X^\lambda), \quad
{h} \mapsto {h} |_{X^\lambda} + \sum_{n,d} (q^d/n!) \ev_{n+1,*}
\ev_1^* {h} \cup \ldots \cup \ev_n^* {h} .$$
\noindent Let $\pi_{G_\lambda}^G: \Lambda_X^{G_\lambda} \to \Lambda_X^G$ be the
canonical map of Novikov rings induced by $H_2^{G_\lambda}(X) \to
H_2^G(X)$. 

\begin{lemma}  Suppose that stable=semistable for $\lambda$-fixed gauged
maps.  Then
\begin{equation} 
 \label{qrestrict} 
\tau_{X,\ti{X}^t,G,\lambda} = \pi_{G_\lambda}^G \circ
\tau_{X^\lambda,\ti{X}^t | X^\lambda,G_\lambda/\C^\times_\lambda}
\circ \iota_\lambda. \end{equation}
\end{lemma}

\begin{proof} 
Decomposing the fixed point locus according to the number of 
markings on each bubble tree gives an isomorphism 
\begin{eqnarray*} 
 \ol{\M}^{G}_n(C,X,\ti{X},\lambda) &\cong& \bigcup_{i_1 + \ldots + i_r
   = n} \prod_{j=1}^r ( \{ \pt \} \cup \ol{\M}_{0,i_j +
   1}(X)^{\C^\times_\lambda}) \times_{(X^\lambda)^r}
 \ol{\M}_r^{G_\lambda,\fr}(C,X^\lambda) / G_\lambda^r \\ &\cong&
 \bigcup_{i_1 + \ldots + i_r = n} \prod_{j=1}^r ( \{ \pt \}/G_\lambda
 \cup \ol{\M}_{0,i_j + 1}(X)^{\C^\times_\lambda}/G_\lambda )
 \times_{(X^\lambda/G_\lambda)^r}
 {\M}_r^{G_\lambda/\C^\times_\lambda}(C,X^\lambda)
     \end{eqnarray*}
where $\{ \pt \}$ represents a trivial bubble tree attached at the
$j$-th node on the principal component.  It follows that integration
over $\ol{\M}^{G}_n(C,X,\ti{X},\lambda)$ is given by push-forward of
$\ev_1^* {h} \cup \ldots \ev_{i_j}^* {h}$ over each
$$\ev_{i_j+1} : \{ \pt \} \cup \ol{\M}_{0,i_j + 1}(X)^{\C^\times_\lambda}/G_\lambda \to X^\lambda/G_\lambda$$
followed by integration over $\ol{\M}_r^{G_\lambda}(C,X^\lambda)$, or
more precisely, $\C^\times_\lambda$-equivariant integration over the
Deligne-Mumford stack
$\ol{\M}_r^{G_\lambda/\C^\times_\lambda}(C,X^\lambda)$ (for which
stable=semistable). 
\end{proof} 

It will be important for our induction argument later that the rank of
the structure group for the fixed point contributions is less than the
rank of the original group.  This allows an approach to results such
as quantum abelianization and quantum Lefschetz by induction on the
rank of the structure group.  More precisely, there exists a canonical
isomorphism
$$\M^{G_\lambda}(C,X^\lambda) \to
\M^{G_\lambda/\C^\times}(C,X^\lambda) .$$
Indeed via the projection map $G_\lambda \to
G_\lambda/\C^\times_\lambda$ any gauged map to $X^\lambda/G$ defines a
gauged map to $X^\lambda/(G_\lambda/\C^\times_\lambda)$ and we obtain
a map
\begin{equation} \label{proj} \M^{G_\lambda}(C,X^\lambda) \to
\M^{G_\lambda/\C^\times}(C,X^\lambda) .\end{equation}
Up to finite cover the exact sequence
$$ 1 \to \C^\times_\lambda \to G_\lambda \to
G_\lambda/\C^\times_\lambda \to 1 $$
splits.  Given a gauged map to $X^\lambda/(G_\lambda/\C^\times))$, let
$c$ denote the weight of the $C^\times_\lambda$-action on $\ti{X}|
X^\lambda$.  Taking the bundle $\C^\times_\lambda$-bundle with first
Chern class $-c$ defines the inverse map to \eqref{proj}.

Finally we remove the requirement that stable=semistable for
linearized gauged maps.  This requirement can be weakened to
stable=semistable for gauged maps with respect to the two
linearizations being compared, using moduli stacks of parabolic gauged
maps:

\begin{proposition} \label{irrat2} 
 If stable=semistable then the moduli stack
 $\ol{\M}_n^G(C,X,\ti{X}_-,\ti{X}_+,\nu)$ of parabolic gauged maps is
 a smooth, proper Deligne-Mumford stack with a perfect relative
 obstruction theory.  For generic parabolic weights $\mu$,
 stable=semistable and so the conclusion holds.
\end{proposition} 

The following proposition is proved in exactly the same way as
Proposition \ref{irrat}.

\begin{corollary}  
If stable=semistable for $\ti{X}_\pm$ then for generic parabolic
weights $\nu$ the stack $\ol{\M}_n^G(C,X,\ti{X}_-,\ti{X}_+,\nu)$ of
Mundet-semistable gauged maps with parabolic structure is a smooth,
proper Deligne-Mumford stack with a perfect relative obstruction
theory containing $G/B$-bundles over $\ol{\M}_n^G(C,X,\ti{X}_\pm)$ as
fixed point components.
\end{corollary}  

In other words, at the cost of adding a parabolic structure we may
always obtain a master space that is a Deligne-Mumford stack with a
perfect relative obstruction theory.

\section{Quantum Witten localization} 
\label{put} 

In this section we combine the area-dependence studied in the previous
section with large and small area limit theorems from \cite{qk3},
\cite{small} to obtain a proof of the quantum Witten localization
formula \eqref{qwitten}.  The gauged potential and the graph potential
of the quotient are related by the {\em adiabatic limit theorem} of
\cite{qk3} (which is a generalization of an earlier result of
Gaio-Salamon \cite{ga:gw}).  

\subsection{Quantum Kirwan map }

In this section we recall the quantum Kirwan map $\kappa_{X,G}$ of
\eqref{qkirwan}.  The map $\kappa_{X,G}$ is defined by virtual
integration over a moduli stack {\em scaled affine gauged maps} to
$X$.

\begin{definition} {\rm (Affine gauged maps)}  Let $n \ge 0$ be an integer.   An $n$-marked affine gauged map is a tuple 
$$( P \to C, u: C \to P(X), \lambda : C \to \P(\omega_C \oplus \C),
  \ul{z} = (z_0,\ldots, z_n) )$$
where $C$ is a twisted balanced curve as in orbifold Gromov-Witten
theory \cite{agv:gw}, $P \to C$ is a principal $G$-bundle, $\omega_C$
is the dualizing sheaf on $C$, and $\lambda$ is a section of its
projectivization $\P(\omega_C \oplus \C)$ which satisfies a certain
{\em monotonicity condition}: on any maximal non-self-crossing path of
components $C_0,C_1,\ldots, C_l$ of $C$ starting with the component
$C_0$ containing $z_0$, $\lambda|C_i$ is non-zero and finite on
exactly one component $C_i$, on which $\lambda$ has a single double
pole.  Such a map is semistable if $u$ takes values in $X \qu G$ on
the locus $\lambda^{-1}(\infty) \subset C$, the bundle $P$ is trivial
on the locus $\lambda^{-1}(0)$, $z_0 \in \lambda^{-1}(\infty)$ while
$z_1,\ldots, z_n \in \lambda^{-1}(< \infty)$ and the datum admits no
automorphisms: each component on which the scaling $\lambda$ is finite
and non-zero resp. zero or infinite and on which $(P,u)$ is
trivializable has at least two resp. three special point.
\end{definition} 

We introduce notation for moduli stacks and evaluation maps.  Each
component $\ol{\M}_{n,1}^G(\C,X,d)$ of homology class $d \in
H_2^G(X,\Z)/\on{torsion}$ and $n$ markings has evaluation maps
$$ \ev_\infty \times \ev: \ol{\M}_{n,1}^G(\C,X,d) \to (X \qu G) \times
(X/G)^n \quad (P,u,\lambda,\ul{z}) \mapsto (u(z_0), \ldots, u(z_n))
.$$
The formula for $\kappa_{X,G}$ is
$$ \kappa_{X,G}({h}) = \sum_{n \ge 0, d} (q^d/n!) \ev_{\infty,*}
\ev^* ({h},\ldots,{h}) .$$

We remark that if the target satisfies an equivariant Fano condition then the
derivative of the quantum Kirwan map at zero is homomorphism of small
quantum cohomologies.  Here equivariantly Fano means that the
equivariant first Chern class $c_1^G(TX)$ of $X$ is the same as the
first Chern class $c_1^G(\ti{X})$ of the linearization.  It follows
that the first Chern class $c_1^G(TX)$ pairs positively with any curve
class for which there is a generically semistable map to the quotient
stack.  In the Fano case, the moduli stacks $\ol{\M}_n^G(\C,X,d)$ have
dimension
$$ \dim \ol{\M}_n^G(\C,X,d) \ge  \dim(X) + 2 .$$  
For ${h} \in H^{\leq 2}_G(X)$, the push-forwards to $H(X \qu G)$ under
$\ev_\infty$ have degree larger than $2 \dim(X \qu G)$ \cite{qk3} for
$d \neq 0$.  Hence
\begin{equation} \label{eqfano}
 \kappa_{X,G}(0) = 0, \quad D_0 \kappa_{X,G}: T_0 QH_G(X) \to T_0 QH(X \qu G) .\end{equation} 
For similar reasons, a lower bound on $m$ the minimal Chern number
$(d,c_1^G(X))$ for classes $d \in H_2^G(X)$ realized by stable
affine gauged maps $u: \P(1,r) \to X/G$ implies that $D_0
\kappa_{X,G}$ has no quantum corrections on classes of degree at most
$2m$.  This ends the remark. 

\begin{example} {\rm (Quantum Kirwan map for the scalar multiplication on affine space)} 
Let $G = \C^\times$ act on $X = \C^k$ by scalar multiplication, so
that $X \qu G = \P^{k-1}$.  We have 
$$T_0 QH_G(X) = \Lambda_X^G[\xi],$$
with \(\xi\) the equivariant parameter, while 
$$T_0 QH(X \qu G) = \Lambda_X^G[\omega]/( \omega^k - q),$$
with \(\omega\in H(\P^{k-1})\) the standard hyperplane class.  By the
previous remark $\kappa_{X,G} (0) = 0$ and
$$ D_0 \kappa_{X,G}(\xi^l) = \omega^l, \quad  l < k .$$ 
A special case of the main result of \cite{surject} (quantum
Stanley-Reisner relations) implies that
$$ D_0 \kappa_{X,G}(\xi^k) = q. $$
Hence $D_0 \kappa_{X,G}$ is surjective and 
$$ T_0 QH(X \qu G) = T_0 QH_G(X) / \ker D_0 \kappa_{X,G} =
\Lambda_X^G[\xi]/ (\xi^k - q) $$
as expected.
\end{example} 

\subsection{Adiabatic limit theorem} 

The following theorem describes the relationship between the gauged
potential and the graph potential of the quotient.  Let $\rho$ be a
positive integer and consider and the family of linearizations
$\ti{X}^\rho$ with $\rho \to \infty$.

\begin{theorem} 
\label{largerel}
{\rm (Adiabatic limit theorem \cite{qk3})} If stable=semistable for
the action of $G$ on $X$ then stable=semistable for gauged maps for
$\rho$ sufficiently large (more precisely, for any class $d \in
H_2^G(X,\Z)$ there exists an $r > 0$ such that $\rho > r$ implies
stable=semistable) and
$$ {\tau}_{X \qu G} \circ {\kappa_{X,G}} = \lim_{\rho \to \infty}
{\tau}_{X,G} : QH_G(X) \to \Lambda_X^G . $$
If $C$ is a genus zero curve equipped with a $\C^\times$-action, then
the same equality holds for $\C^\times$-equivariant potentials
$\tau_{X \qu G}^{\C^\times}, \tau_{X,G}^{\C^\times}$.  Similarly for
the localized gauged potentials, consider the $\C^\times$-equivariant
version of the quantum Kirwan map
$$ \kappa_{X,G}: QH_G(X) \to QH_{\C^\times}(X \qu G) $$
Then \cite[Theorem 1.6]{qk1}
\begin{equation} \label{locad}
 \tau_{X \qu G,\pm} \circ \kappa_{X,G} = \tau_{X,G,\pm}: QH_G(X) \to
 QH_{\C^\times}(X \qu G). \end{equation}
\end{theorem} 

\noindent In other words, the diagram
\begin{equation} \label{classdiag}
  \begin{diagram} \node{QH_G(X)} \arrow{se,b}{\tau_{X,G}}
  \arrow[2]{e,t}{{\kappa}_{X,G}} \node{} \node{QH( X \qu G)}
  \arrow{sw,b}{\tau_{X \qu G}} \\ \node{} \node{\Lambda_X^G} \node{}
   \end{diagram}\end{equation}
commutes in the limit $\rho \to \infty$.  

\begin{remark} \label{allows} In semi-Fano cases the localized adiabatic limit equation \eqref{locad}
allows to solve for the quantum Kirwan map on divisor classes.  The
constant term in $\tau_{X,G,\pm}$ (with respect to the equivariant
parameter) arises from configurations with $\phi_\pm$, $d_\pm$ and
$n_\pm$ vanishing, which are constant maps to $X \qu G$ with no
markings.  It follows that the expansion of $\tau_{X,G,\pm}$ in powers
of the inverted equivariant parameter $1/\zeta$ is
$$ \tau_{X,G,\pm} = 1 + \tau_{X,G,\pm}^{(1)}/\zeta +
\tau_{X,G,\pm}^{(2)}/\zeta^2 + \tau_{X,G,\pm}^{(3)}/\zeta^3 + \ldots
.$$
Similarly the constant contribution to $\tau_{X \qu G,\pm}$ is $1$,
corresponding to configurations with no markings at $0$ or $\infty$,
while the contribution from a single marking again involves no
bubbles, with a single marking at $0$ or $\infty$ with Euler class
$\zeta$ and so is the $\on{Id}/\zeta$.  It follows that the expansion
of $\tau_{X \qu G, \pm}$ in powers of $\zeta^{-1}$ is
$$ \tau_{X \qu G,\pm} = 1 + \on{Id}/\zeta 
+ \tau_{X \qu G,\pm}^{(2)}
/\zeta^{2}
 + \tau_{X \qu G,\pm}^{(3)}/ \zeta^{3} + \ldots : QH_G(X)
\to QH_{\C^\times}(X \qu G) .$$
Similarly $\kappa_{X,G}$ admits an expansion in powers of $\zeta$,
$$ \kappa_{X,G} = \kappa_{X,G}^{(0)} + \kappa_{X,G}^{(1)} \zeta +
\kappa_{X,G}^{(2)} \zeta^2 + \ldots .$$
Suppose that $X$ is equivariantly semi-Fano in the sense that
$c_1^G(X)$ is non-negative on the class of every affine gauged map.  Then
for reasons of degree
$$ \kappa_{X,G}( QH^{\leq 2}_G(X) ) \subset QH^{\leq 2}(X \qu G), \quad
\kappa_{X,G}^{(\ge 2)} (QH^{\leq 2}_G(X)) = 0 .$$
The constant term in the relation
$$ (1 + \on{Id}/\zeta  + \tau_{X \qu G,\pm}^{(2)} /\zeta^{2})
\circ(\kappa_{X,G}^{(0)} + \kappa_{X,G}^{(1)} \zeta)+ O(\zeta^{-2}) =
1 + \tau_{X,G,\pm}^{(1)}/\zeta + O(\zeta^{-2}) \quad \text{on}
\ QH^{\leq 2}(X \qu G)$$
implies that $\kappa_{X,G}^{(1)}$ vanishes on $QH^{\leq 2}_G(X)$.
Hence
\begin{equation} \label{same} \kappa_{X,G}^{(0)} = 
\tau_{X,G,\pm}^{(1)} \quad \text{on} \ \  QH^{\leq 2}_G(X) .\end{equation}
This ends the remark.
\end{remark}

\label{smallarea}
The limit of the Mundet semistability condition in which the
linearization goes to zero is studied in the paper \cite{small}.  In
this limit, the bundle must be semistable and so the moduli stack of
gauged maps is a quotient of the moduli space of parametrized stable
maps to $X$.  Theorem \ref{qwform} now follows from Theorems
\ref{gwall} and Theorem \ref{largerel}.

\subsection{Localization for convex varieties} 

A slightly modified version of the quantum Witten localization formula
holds in quasiprojective cases under a {\em convexity} assumption.

\begin{definition} A finite dimensional complex
$G$-vector space $V$ will be called {\em convex} if there exists a
  central one-parameter subgroup $\phi_\lambda: \C^\times \to G$ such
  that $X$ has positive weights for the induced action of $\phi$,
$$ V = \bigoplus_\mu V_\mu, \quad (\mu,\lambda) > 0 .$$
Given a convex $G$-vector space, the {\em projectivization} of $V$ is
the quotient
$$ \ol{V} = ((V \times \C)^\times - \{ (0 , 0) \})/\C^\times $$
where $\C^\times$ acts on $\C$ with weight one.  Thus $\ol{V}$ is a
weighted projective space and contains $V$ as an open subset.  A
quasiprojective $G$-variety $X$ is {\em convex} if there exists a
projective morphism $\pi: X \to V$ to a convex $G$-vector space $V$.
\end{definition} 

The following is a simple application of the technique called {\em
  symplectic cutting} in the literature \cite{le:sy2}:

\begin{lemma}   Any convex $G$-variety $X$ admits a $G$-equivariant 
compactification $\ol{X}$ by adding single $\C^\times_\lambda$-fixed
divisor.
\end{lemma} 

\begin{proof} 
Let $\ti{X} \to X$ denote the given linearization on $X$ and
$\ti{X}(k)$ the linearization on $X \times \C$ obtained by twisting by
the $\C^\times$-character with weight $k$.  Consider the git quotient
$$ \ol{X} = (X \times \C) \qu \C^\times .$$
The inverse image of $ (0,0) \in V \times \C$ is unstable, for
sufficiently large $d$.  Thus the proper morphism $X \to V$ induces a
proper morphism $\ol{X}$ to $\ol{V}$.  In particular, the quotient
$\ol{X}$ is also proper.  The $G$ action on $X \times \C$ given by
$g(x,z) = (gx,z)$ descends to a $G$-action on $\ol{X}$, and restricts
to the given action on the open subset $X \subset \ol{X}$.
\end{proof} 

\begin{corollary}  \label{disjointcor} Let $d \in H^2_G(\ol{X})$ be a class
that pairs trivially with the divisor class $[\ol{X} - X] \in
H_2^G(\ol{X})$.  Then for $k \gg 0$ the moduli stack
$\ol{\M}_n^G(C,\ol{X},d)$ consists of maps whose images are disjoint
from $(\ol{X} - X)/G$.  Similarly, if $\ti{X}_\pm \to X $ are two
different linearizations then for $k \gg 0 $ the moduli stack
$\ol{\M}_n^G(C,\ol{X},d)$ consists of maps whose images are disjoint
from $(\ol{X} - X)/G$.
\end{corollary}

\begin{proof} The intersection number of any curve
  $u: \P^1 \to \ol{V}$ contained in $\ol{V} - V$ with $\ol{V} - V$ is
  non-negative.  Indeed $\ol{V} - V$ has ample normal bundle in
  $\ol{V}$ being a prime invariant divisor in a weighted projective
  space.  On the other hand, there are no stable gauged maps
  $C \to X/G$ with image in $(\ol{V} - V)/G$ for sufficiently large
  $d$, Indeed the trivial reduction $\sigma$ together with the
  generator $\lambda$ of the one-parameter subgroup $\C^\times$ has
  weight $\mu(\sigma,\lambda) \to \infty $ as $d \to \infty$.
  Combining these observations let $v: \hat{C} \to \ol{V}/G$ be a
  stable gauged map intersecting $(\ol{V} - V)/G$.  Then the
  intersection number $\# u^{-1}(P(\ol{V}- V)) > 0 $ is positive and
  equal to the pairing $(d, [\ol{V} - V]) \in \Q$ of
  $d \in H_2^G(X,\Q)$ with $[\ol{V} - V] \in H^2_G(\ol{V},\Q)$.  The
  latter vanishes by assumption, a contradiction.
\end{proof} 

The corollary implies that the wall-crossing formula also holds for
convex varieties by applying the formula to the compactified variety
with compactifying divisor sufficiently far away at infinity. However,
the quantum Witten localization formula does not hold because,
eventually, the compactifying divisor will make a contribution in the
localization formula.  The following alternative argument gives a
formula similar to that in quantum Witten localization.  Let $\chi$ be
a character of $G$ that is negative on the one-parameter subgroup
generated by $\xi$ and $\ul{\C}_\chi$ the corresponding trivial line
bundle over $X$.  Consider the piecewise linear path of linearizations
$\ti{X}_\rho \to X$ obtained by shifting by multiples of the character
$\chi$:
\begin{equation}
 \label{convexpath}
\ti{X}_\rho = \begin{cases} \ti{X} \otimes \ul{\C}_\chi^{\rho^{-1} -
    1} & \rho \leq 1 \\ \ti{X}^\rho & \rho \geq 1 \end{cases} .\end{equation}

\begin{lemma}  \label{emptylem}
For any homology class $d \in H_2^G(X,\Z)$, the moduli stack
$\ol{\M}_n^G(C,X,\ti{X}_\rho,d)$ is empty for $\rho \gg 0$.
\end{lemma}

\begin{proof}  Let $\sigma: C \to P/G$ be the trivial parabolic reduction, and $\lambda$
the generator of the one-parameter subgroup in the definition of
convexity.  Given a gauged map $v: \hat{C} \to X/G$, the associated
graded pair $\Gr(P),\Gr(u) $ for $(\sigma,\lambda)$ projects to the
origin in $V/G$.  The Mundet weight picks up a term $(\rho^{-1} -
1)(\chi,\lambda)$ which goes to infinity as $\rho \to 0$.  Hence there
are no Mundet-semistable gauged maps with class $d$, for $\rho$
sufficiently small.
\end{proof} 

\begin{theorem}
 \label{vspace}  {\rm (Quantum localization for convex
varieties)} Let $X$ be a convex $G$-variety, $C$ a genus zero curve,
 and suppose that stable=semistable for the $G$-action on $X$, for
 gauged maps with linearization $\ti{X}$, and for polarized gauged maps
 for the path $\ti{X}_\rho$.  Then
 \begin{equation} \label{qwitten3} \tau_{X \qu G} \circ \kappa_{X,G} =
   \sum_{ [\lambda] \neq 0,\rho } \frac{|W_{\lambda}|}{|W_{\C
       \lambda}|}  \tau_{X,G,\lambda,\rho} .\end{equation}
where the sum is over equivalence classes $[\lambda]$ of
unparametrized one-parameter subgroups generated by $\lambda
\in \g$.    Similarly for the localized graph potentials 
from \eqref{infty}
$$ \tau_{X \qu G,\pm}^\infty \circ (\kappa_{X,G} \times \kappa_{X \qu
  G}) = \sum_{ [\lambda] \neq 0,\rho } 
\frac{|W_{\lambda}|}{|W_{\C
    \lambda}|}  \tau_{X,G,\lambda,\pm}^\rho.$$
\end{theorem}

\begin{proof} This is a combination of the adiabatic limit theorem
  \ref{largrel}, the wall-crossing formula Theorem \ref{gwall}, the
  vanishing of the invariants for large $\rho$ in Lemma
  \ref{emptylem}.  The application of these results to the non-proper
  variety $X$ is justified by the relationship between invariants of
  the compactification $\ol{X}$ with those of $X$ in Corollary
  \ref{disjointcor}.
\end{proof}

\begin{example} \label{affine}  
{\rm (Quantum Witten localization for the scalar multiplication on
  affine space)} To explain the notation we use \eqref{qwitten3} to
compute the three-point Gromov-Witten invariants of projective space
using quantum Witten localization.  Suppose that $G = \C^\times$ acts
diagonally on $X = \C^k$ so that
$$ X \qu G = \C^k \qu \C^\times = \P^{k-1} .$$ 
 We have 
$$H_2^G(X,\Z) \cong H_2(X \qu G) = \Z [\P^1], \quad H^2_G(X,\Z) \cong
 H^2(X \qu G) = \Z \omega $$
 where $\omega$ is the hyperplane class, the image of the equivariant
 generator $\xi \in H^2_G(X,\Z)$ under the Kirwan map.  We compute the
 class $d = 1$ three-point invariants using quantum Witten
 localization.  Let $\beta \in H^6(\ol{\M}_3(C))$ be the fundamental
 class, whose insertion fixes the positions of the marked points in
 \eqref{ins}.  We identify $QH_G(X) \cong \Lambda_X^G[\xi]$.  Consider
 the three-point invariants with insertions
 $\xi^a, \xi^b, \xi^c \in QH_G(X) \cong S(\g)^\dual .$ Since
 $c_1^G(X)$ is at least $2k$ on classes $d > 0$, the derivative
 $D_0 \kappa_{X,G}$ of the quantum Kirwan map has no quantum
 corrections by \eqref{eqfano}.  The image of $\xi^a,\xi^b,\xi^c$
 under $D_0 \kappa_{X,G}$ is equal to $ \omega^a, \omega^b, \omega^c$
 respectively.  We consider a path $\ti{X}_\rho$ obtained by shifting
 by a negative character $\chi$; this means that in the fixed point
 formula we take the residue with respect to $-\xi$, see \cite{wall}.
 By the formula \eqref{qwitten3},
\begin{eqnarray*} 
\sum_{d \ge 0} q^d \lan \omega^a, \omega^b, \omega^c \ran_{0,d} &=&
\tau_{X \qu G}^3(\omega^a,\omega^b,\omega^c,\beta) \\&=& -
\sum_{\rho,[\lambda]}
\tau^3_{X,\ti{X}_\rho,G,\lambda}(\xi^a,\xi^b,\xi^c,\beta). \end{eqnarray*}
There is a unique $G$-fixed point in $X$.  The $G$-bundle $P$ with
first Chern class $d = 1$ together with the zero section $u \in
H^0(C,P \times_G X)$ forms a Mundet semistable map for a unique value
of the parameter $\rho$.  For $d = 1$ the index bundle and its Euler
class are
$$ \Ind(T(X/G)) = H^0(\mO(k)^{\times} \times_{\C^\times} \C^k) \cong
\C^{2k} , \quad \eps_+(\Ind(T(X/G))) = \xi^{2k} .$$
The unique fixed point contribution
\begin{eqnarray*} 
\tau^3_{X,\ti{X}_\rho,G,\lambda}(\xi^a,\xi^b,\xi^c,\beta) 
&=& q \Resid_{-\xi}   \frac{\xi^{a + b + c}}{\xi^{2k}}  \\
&=& \begin{cases} q &  a + b + c = 2k-1 \\ 0 &  \text{otherwise} .\end{cases}  
\end{eqnarray*}
We obtain
$$\lan \omega^a, \omega^b, \omega^c \ran_{0,1} =
\begin{cases} 1 & a + b +c   = 2k - 1,  \\ 
              0 & \text{otherwise} \end{cases} $$
as expected.    
\end{example} 

\section{Abelianization for Gromov-Witten invariants}  
 
In this section we prove Theorem \ref{qabel} from the introduction.
Abelianization is first proved for graph potentials, then deduced for
qde solutions.  The section ends with the examples of the moduli
spaces of odd numbers of points on a projective line and the
Grassmannian.

\subsection{Abelianization for graph potentials}

\begin{proof}[Proof of Theorems \ref{qabel} and \ref{abelloc}]
We first prove the Theorem in the case stable=semistable for
linearized gauged maps.  We take as the inductive hypothesis that
Theorem \ref{qabel} holds for any group of rank less than $\dim(G)$.
We wish to compare the fixed point contributions in the quantum Witten
localization formulas
\begin{equation} \label{qwitten1} \tau_X^G - \tau_{X \qu G} \circ \kappa_{X,G} = \sum_{
  [\lambda] \neq 0,\rho } \tau_{X,\ti{X}^\rho,G,\lambda} \end{equation}
and 
\begin{equation}  \label{qwitten2}  \tau_X^{T,\/g/\t} - \tau_{X \qu T}^{\g/\t}  \circ \kappa^{\g/\t}_{X,T} = \sum_{
  [\lambda] \neq 0 ,\rho} \tau^{\g/\t}_{X,\ti{X}^\rho,T,\lambda} .\end{equation}
In the version for $T$, both the traces and quantum Kirwan maps have
been twisted by the Euler class of the index of $\g/\t$.  Now
$\tau_X^{T,\g/\t}$ resp. $\tau_X^G$ is defined by integration over
$\ol{\M}_n(C,X) \qu T$ resp. $\ol{\M}_n(C,X) \qu G$.  This is
essentially the setting considered by Martin \cite{mar:sy}.  In
Gonz\'alez-Woodward \cite[Chapter 5]{small} we show
$$ \tau_X^G = |W|^{-1} \pi_T^G \circ  \tau_X^{T,\g/\t} $$
either by Martin's argument, if the moduli spaces of stable maps are
smooth and the virtual fundamental classes are the usual ones, or by a
virtual version of Martin's argument if the moduli spaces of stable
maps are only virtually smooth.  We note that the virtual non-abelian
localization formula used in \cite{small} had a gap in the proof,
which was fixed by Halpern-Leistner \cite{hl:qs}.

We review the argument for abelianization in the small-area limit
briefly from \cite{small}.  Write $\chi_G$ for the Euler
characteristic on the stack $S/G$ and $\chi_T$ for the Euler
characteristic on the stack $S/T$.  By the Weyl character formula
$$ \chi_G(F) = |W|^{-1} \chi_T(F \otimes \ul{\on{Alt}}(\g/\t)) .$$
Here $\ul{\on{Alt}}(\g/\t)$, the trivial sheaf with values in the
exterior algebra $\on{Alt}(\g/\t)$, is the K-theory Euler class given
by the trivial bundle with fiber the $T$-representation
$\on{Alt}(\g/\t)$ with character $\prod_{\alpha \in \RR} (1 -
t^{\alpha})$.  For sufficiently positive bundles we have
\begin{equation} \label{chiab}
 \chi(S \qu G, F \qu G) = |W|^{-1} \chi_{T}(S \qu T, F \otimes
 \ul{\on{Alt}}(\g/\t) \qu T)
 \end{equation}
Using virtual Riemann-Roch in Tonita \cite{to:rr} and taking limits as
we obtain
$$ \int_{S \qu G} \Ch(F \qu G) = |W|^{-1} \int_{S \qu T} \Ch(F \qu T)
\cup \Eul(\g/\t)) $$
where $\Ch(F) \in H(S)$ is the Chern character and $\Eul(\g/\t)$ the
Euler class in cohomology. This implies the abelianization formula in
the small-area limit for Chern characters:
$$ \int_{S \qu G} \kappa_{X,G}( h) = |W|^{-1} \int_{S \qu T}
\kappa_{X,T}( h \cup \Eul(\g/\t)) , \quad \forall h \in
\HH_G(X).$$

Using abelianization in the small-area limit to prove abelianization
it suffices to show abelianization for the right-hand-sides in
\eqref{qwitten1}, \eqref{qwitten2}.  Each fixed point component
$X^\lambda$ for the $G$-action corresponds to $|W/W_{\C\lambda}|$
fixed point components $X^{w \lambda}, w \in W/W_{\C\lambda}$ for the
$T$-action.  The identity we wish to show is
\begin{equation} \label{want}  
  \tau_{X,\ti{X}^\rho,G,\lambda} = |W_\lambda|^{-1} \pi_T^{G_\lambda} \circ
  \tau^{\g/\t}_{X,\ti{X}^\rho,T,\lambda} .
\end{equation} 
In the case $G_\lambda$ is abelian the group $W_\lambda$ is trivial
and so the equality holds automatically.  More generally the equation
\eqref{qrestrict} gives
$$ \tau_{X,\ti{X}^t,G,\lambda} = \tau_{X^\lambda,G_\lambda/\C^\times_\lambda} \circ
\iota_\lambda $$
and 
$$ |W_\lambda|^{-1} \tau_{X,\ti{X}^t,T,\lambda} = 
|W_\lambda|^{-1}
\tau_{X^\lambda,T/\C^\times_\lambda} 
\circ \iota_\lambda . $$
By the inductive hypothesis, 
$$ \tau_{X^\lambda,G_\lambda/\C^\times_\lambda} = |W_\lambda|^{-1} \pi_T^{G_\lambda}
\circ \tau^{\g_\lambda/\t}_{X^\lambda,T/\C^\times_\lambda} .$$
Equation \eqref{want} follows.  See Guillemin-Kalkman \cite[Section
  4]{gu:ne} for similar arguments involving recursive applications of
fixed point formulae.  The equality with $\tau_{X \qu T}^{\g/\t} \circ
\kappa_{X,T}^{\g/\t} \circ \Restr_T^G$ follows from by combining
wall-crossing with Theorem \ref{largerel}.  Theorem \ref{abelloc} is
proved by the same argument applied to the $\C^\times$-fixed locus in
$\ol{\M}^G(C,X,\ti{X}_-,\ti{X}_+)$, where $\C^\times$ acts on $C \cong
P^1$ by the action with weights $-1,1$.

To prove the Theorems \ref{qabel} and \ref{abelloc} in the general
case when semistable $\neq$ stable, we find a master space for which
stable=semistable by adding a parabolic structure as in Proposition
\ref{irrat2}.  Given stability parameters $\rho_\pm$ and corresponding
linearizations $\ti{X}_\pm = \ti{X}^{\rho_\pm}$, let
$\ol{\M}_n^G(C,X,\ti{X}_-,\ti{X}_+,\nu)$ of polarized gauged maps with
parabolic structure.  As in Proposition \ref{irrat}, for generic $\mu$
stable=semistable for $\ol{\M}_n^G(C,X,\ti{X}_-,\ti{X}_+,\nu)$.  So
$\ol{\M}_n^G(C,X,\ti{X}_-,\ti{X}_+,\nu)$ is also a proper
Deligne-Mumford stack with perfect relative obstruction theory.
Localization on the stack $\ol{\M}_n^G(C,X,\ti{X}_-,\ti{X}_+,\nu)$
produces a wall-crossing formula, whose fixed point contributions are
the $\C^\times$-fixed components in
$\ol{\M}_n^G(C,X,\ti{X}_-,\ti{X}_+,\nu)$.  By induction on the rank of
$G$, we obtain the identity
\begin{multline} 
 \int_{[\ol{\M}_n^G(C,X,\ti{X}^\rho,\nu,d)]} \ev^* {h} \cup \pi^* \Eul(L_\pi) \\ =
 \int_{[\ol{\M}_n^T(C,X,\ti{X}^\rho,\nu,d)]} \ev^* {h}  \cup \eps(\g/\t) \cup \pi^* \Eul(L_\pi)/ |W| \end{multline}
for any $\rho$ for which stable=semistable for gauged maps.  Hence
using \eqref{fiberint},
$$ \int_{[\ol{\M}_n^G(C,X,\ti{X}^\rho,d)]} \ev^* {h} =
\int_{[\ol{\M}_n^T(C,X,\ti{X}^\rho,d)]} \ev^* {h} \cup \eps(\g/\t) /
|W| $$
as claimed.
\end{proof}

\subsection{Abelianization for convex varieties} 

\label{vspace2}  
 Continuing Theorem \ref{vspace}, we extend the abelianization results
 to convex $G$-varieties.  Replacing the quantum Witten localization
 formula \eqref{qwitten} by the alternate formula \eqref{qwitten3}
 obtained from wall-crossing by shifting by a character, the same
 arguments go through and imply the formulas \ref{qabel} and
 \ref{abelloc}.  This also gives alternative argument for
 abelianization for git quotients $X \qu G$ of projective $X$ by $G$
 in the case that $G$ has a non-trivial center, but not in the case
 that $G$ is simple and non-abelian. This allows us, finally, to give
 some applications.

\begin{example} \label{grass} 
{\rm (Grassmannians)} In this example we reproduce the results on
Grassmannians from Bertram et al \cite{be:qu}.  For positive integers
$r < k$ let $X = \Hom(\C^r,\C^k)$ be the space of linear maps from
$\C^r$ to $\C^k$. Let $G = GL(r)$ act on $X$ by composition.  For the
linearization $\lambda = c_1^G(X) \in H^2_G(X) \cong \Q$, the
semistable maps are those with full rank and so
\begin{eqnarray*} 
X \qu G &=& \{ x \in \Hom(\C^r,\C^k) | \on{rank}(x) = k \}/ GL(r) \\ 
&=& \Gr(r,k) \end{eqnarray*}
is the Grassmannian of $r$-dimensional subspaces in $\C^k$ with
$H_2^G(X,\Z) \cong \Z$.  There exist stable maps to $X \qu G$ of class
$d$ only if $d \ge 0$. Although $X$ is not compact, the weights for
the central $\C^\times$ action are all one, and so the git quotients
are proper.  Since $X$ is equivariantly Fano, the map $\kappa_{X,G}$
is trivial on $H^2_G(X)$ and may be ignored by \eqref{eqfano}.  By
abelianization the localized graph potential $\tau_{X \qu G,\pm}$ for
the Grassmannian has restriction to $H^2(X \qu G)$ given by
\begin{equation}\label{hv1} 
 \tau_{X \qu G,\pm} = \mu_T^G \circ \tau_{X,T,\pm} \circ
 \Restr_T^G \end{equation}
Let $\theta_1,\ldots, \theta_r$ denote the standard basis of
$H_T^2(X,\Z) \cong \Z^r$.  The localized gauged potential
$\tau_{X,T,\pm}$ is Givental's $I$-function and by \cite{gi:eq} given by
\begin{equation} \label{hv2}
\tau_{X,T,\pm}^{\g/\t}(t_0 +t_1 \theta_1 + \ldots + t_r \theta_r) = \sum_{\ul{d} \ge 0 } q^{\ul{d}} 
e^{t_0 + (t_1 (\theta_1 + d_1 \zeta) + \ldots + t_{2k+2} (
  \theta_{2k+2} + d_{2k+2}\zeta)) / \zeta}
\tau_{X,T,\pm}^{\g/\t}(\ul{d}) \end{equation}
where
\begin{multline} \label{hv3} \tau_{X,T,\pm}^{\g/\t}(\ul{d}) = 
 \frac{ (-1)^{(k-1)d} \prod_{i < j} ((\theta_i - \theta_j) + (d_i -
   d_j) \zeta) } {\prod_{i < j} (\theta_i - \theta_j) \prod_{i=1}^r
   \prod_{l=1}^{d_i} (\theta_i +l \zeta)^n} \\ = \prod_{i \neq j}
 \frac{ \prod_{ l \leq d_i - d_j} ((\theta_i - \theta_j) + l \zeta) }
      { \prod_{ l \leq 0} ((\theta_i - \theta_j) + l \zeta )}
      \frac{1}{ \prod_{i=1}^r \prod_{l=1}^{d_i} (\theta_i +l
        \zeta)^n}.
\end{multline}
The formula obtained by combining \eqref{hv1}, \eqref{hv2},
\eqref{hv3} was conjectured in Hori-Vafa \cite[Appendix]{ho:mi}.  The
formula was proved in Bertram et al \cite{be:qu} , \cite{be:tw} by
different methods.
\end{example} 

\subsection{Abelianization for quantum products} 

In this section we prove the Theorem \ref{compare} relating the
quantum products for the abelian and non-abelian products, using the
abelianization formula for qde solutions from the previous section.
As an example, we describe the small quantum cohomology ring of the
Grassmannian following Schm\"ashcke \cite{schm:thes}.

We introduce the following differential operator: For any root
$\alpha$ let
$$\alpha_0 = \kappa_{X,T}^{q =0}(\alpha) = c_1(\C_\alpha \qu T) \in
H^2(X \qu T)$$
denote the corresponding class on the symplectic quotient.  Define
$$ \D^{\g/\t} =\prod_{\alpha \in \RR_+} \frac{
  \partial_{\alpha}}{\alpha_0} : \End(\Map(H_T(X), H_{\C^\times \times
  \C^\times}(X \qu T))) ;$$
here the invertibility of $\alpha_0$ may be achieved by adding an
additional equivariant parameter, as in the definition of inverted
Euler classes.  In the equivariant cohomology $H_{\C^\times \times
  \C^\times}(X \qu T)$ the first factor corresponds to the rotation
action on $\C^\times$, while the second factor is that used to define
the Euler-twisted potential in \eqref{epsE}.

\begin{proposition} \label{diffop}   
The $\g/\t$-twisted potential $\tau_{X,T,\pm}^{\g/\t}$ is related to the
untwisted potential by
$$ \tau_{X,T,\pm}^{\g/\t} = \D^{\g/\t} \tau_{X,T,\pm} : QH_T(X) \to
QH_{\C^\times \times \C^\times} (X \qu T) .$$
\end{proposition} 

\begin{proof} 
The addition of $\eps(\g/\t) = \Eul(\Ind(\g/\t))$ into the definition
of $\tau_{X,T,\pm}$ causes the appearance of the additional factors
involving $\alpha_0$:
\begin{eqnarray*} 
 \tau_{X,T,\pm}^{\g/\t} &=& \sum_{d \in H_2^T(X)} q^d \prod_{\alpha > 0}
 \Delta_{d}(\alpha_0,0) \Delta_{d}(-\alpha_0, 0)
 \tau_{X,T,\pm,d} \\ 
&=&  \sum_{d \in H_2^T(X)} q^d \prod_{\alpha > 0} \frac{ \prod_{l =
     -\infty}^{\alpha \cdot d} (\alpha_0 + l \zeta)} { \prod_{l =
     -\infty}^{0} (\alpha_0 + l \zeta)} \frac{ \prod_{l =
     -\infty}^{-\alpha \cdot d} (\alpha_0 + l \zeta)} { \prod_{l =
     -\infty}^{0} (-\alpha_0 + l \zeta)} \tau_{X,T,\pm,d}
 \\ &=& \sum_{d
   \in H_2^T(X)} q^d \prod_{\alpha > 0} \frac{\alpha_0 + (\alpha_0 \cdot
   d) \zeta }{\alpha_0 } \tau_{X,T,\pm,d} \\ &=& \sum_{d \in H_2^T(X)}
 q^d 
\prod_{\alpha \in
   \RR_+} \frac{ \partial_\alpha }{\alpha_0}
\tau_{X,T,\pm,d} 
 = \D^{\g/\t}
 \tau_{X,T,\pm} \end{eqnarray*}
as claimed.
\end{proof} 

\begin{proof}[Proof of Theorem \ref{compare}] 
Let $\Box$ be a $W$-invariant constant-coefficient differential operator
on $QH_2^T(X)$ with symbol $\sigma(\Box)$.  Since $\tau_{X,G,\pm}$ is a
fundamental solution,
\begin{eqnarray*}
 (\Box \tau_{X,G,\pm})({h}) = 0 &\iff& 
 (\Box \tau_{X,T,\pm}^{\g/\t})({h}) = 0\\
 &\iff& D_{h} \kappa_{X,G}
(\sigma(\Box) ) = 0 \in T_{\kappa_{X,G}({h})} QH(X \qu G) .\end{eqnarray*}
On the other hand,
\begin{eqnarray*} 
 \Box \tau_{X,T,\pm}^{\g/\t} = 0 & \iff & \Box \D^{\g/\t} \tau_{X,T,\pm} = 0
 \\ &\iff & D_{h} \kappa_{X,T}( \sigma(\Box \D^{\g/\t})) = 0 \\ &\iff &
 D_{h} \kappa_{X,T}(\sigma(\Box)) D_{h} \kappa_{X,T}(\sigma(\D^{\g/\t}))
 =0 \\ &\iff & D_{h} \kappa_{X,T}(\sigma(\Box)) \in \on{ann}
 D_h \kappa_{X,T}(e_\pm)
\end{eqnarray*} 
as claimed.
\end{proof} 

\begin{example} {\rm (Grassmannians)}  
The relations in $QH( (\P^{n-1})^k)$
are $H_i^n = q, i = 1,\ldots, k$, where $H_i$ is the hyperplane class
on the $i$-th factor.  The Weyl group $W = S_k$ is the $k$-symmetric
group.  The $W$-invariant part of the cohomology ring $QH(
(\P^{n-1})^k)$ is generated by the Schur polynomials
$$ \chi_{\lambda^\dual}(H_1,\ldots,H_n) = \prod_{w \in W} \frac{(-1)^{l(w)}
  H^{w(\lambda + \rho) - \rho}}{ \prod_{i < j} (H_j - H_i) } $$
where $ H^\lambda = H_1^{\lambda_1} \ldots H_k^{\lambda_k}, \quad \rho
= (1,\ldots, k) .$ By the Fano condition $\kappa_{X,T}$ has no quantum
corrections and so
$$ D_0 \kappa_{X,T}(e_\pm) = \pm \prod_{i < j} (H_j - H_i) .$$
Hence for any $\mu \in \Z^k$, we have 
$$ \chi_{\lambda + n \mu} (H_1,\ldots,H_n) -
\chi_{\lambda^\dual}(H_1,\ldots,H_n) \in \on{ann}(\kappa_{X,T}(e_\pm))
\subset H(X \qu T).$$
This implies the relations
$$ \chi_{\lambda + n \mu} = q^{\mu_1 + \ldots + \mu_k} \chi_{\lambda^\dual}
\in QH(\Gr(k,n)), \quad \lambda,\mu \in \Z^k .$$
The leading order terms in these relations are 
$$ \chi_{\lambda + n \mu} = 0, \quad \mu \ge 0 .$$
These are the usual relations in the cohomology of the Grassmannian
describing the cohomology as a truncation of the polynomial
representation ring of $GL(k)$ obtained by setting the Schur
polynomials corresponding to Young diagrams not fitting in the $k
\times (n-k)$ box to zero:
$$ H(G(k,n)) = \Rep(GL(k))/ \lan \chi_{\lambda + n \mu}, \quad \mu > 0
\ran .$$
This implies that the relations above generate the ideal of relations
in the small quantum cohomology ring.  So one obtains the standard
presentation of the quantum cohomology as the Verlinde algebra as in
Bertram-Ciocan-Fontanine-Fulton \cite{be:qm}:
$$ QH(G(k,n)) = \Rep(GL(k))/ (\chi_{\lambda + n \mu} - q^{\mu_1 +
  \ldots + \mu_k} \chi_{\lambda} , \ \ \ \lambda, \mu \in \Z_{\ge 0}^k) .$$
\end{example}  

\section{Quantum Lefschetz for holomorphic symplectic quotients} 

In this section we prove a formula for the qde solution for
hypersurfaces for bundles associated to the semistable locus,
including hypersurfaces defined by the zero level sets of holomorphic
moment maps.  in which the standard techniques introduced for complete
intersections defined by concavex bundles break down.  Suppose that
$X$ is a smooth variety equipped with the action and an equivariant
map
$$ G \times X \to X, \quad \Phi: X \to V .$$
In our examples $\Phi$ will be a holomorphic moment map.  Let
$\ti{X} \to X$ be a linearization, that is, an ample $G$-line bundle.
We denote by $Z = \Phi^{-1}(0)$ the zero level set.  In this section
we give a formula for the graph potential of the git quotient
$Z \qu G$.  Over $X$ we have a natural bundle $X \times V \to X$; we
wish to compare the potentials
$$ \tau_Z \circ r_{Z,G} : QH_G(X) \to \Lambda_X^G, \quad \tau_X^{V}
: QH_G(X) \to \Lambda_X^G $$
where $\tau_X^V$ denotes the potential twisted by the index of the
Euler class of $V$ and $r_{Z,G}$ as before pull-back.

\begin{proof}[Proof of Theorem \ref{qlef}]
We apply the wall-crossing formula in Theorem \ref{gwall} to both
sides: Consider the family of linearizations $\ti{X} \otimes \C_{\rho
  t}$, where $\rho$ is the weight corresponding to the one-parameter
subgroup.  The wall-crossing formulas give a relationship between the
potentials for $\ti{X}_\pm = \ti{X}_{t_\pm}$ at $t = t_\pm$
\begin{equation} \label{wc1} \tau^V_{X,G,1} - \tau^V_{X,G,-1} =
  \sum_{[\lambda],t \in (-1,1)} 
\frac{|W_{\lambda}|}{|W_{\C
    \lambda}|} 
  \tau^V_{X,G,\lambda,t} \end{equation}
\begin{equation} \label{wc2} \tau_{Z,G,1} - \tau_{Z,G,-1} =
  \sum_{[\lambda],t \in (-1,1)} 
\frac{|W_{\lambda}|}{|W_{\C
    \lambda}|} 
  \tau_{Z,G,\lambda,t} \end{equation}%
As in Proposition \ref{irrat2}, at the cost of adding a parabolic
structure we may assume stable=semistable for linearized gauged maps.
By induction on the dimension of the group and \eqref{qrestrict} we
may assume
\begin{equation} \label{walleq} \tau_{Z,G,\lambda,t} = \tau^V_{X,G,\lambda,t}, \quad \forall t \in
(-1,1) .\end{equation}
In the small-area chamber the equality $ \tau_Z^G = \tau_{X}^{G,V} $
holds by standard properties of the Euler class: since $V$ is a
$G$-representation its index bundle
$$ \Ind(V)_{[u:C \to X]} := H^0( u^* TV) \cong V .$$
The locus $ \ol{\M}_n(C,Z,d)$ is the zero set of the map
$\ol{\M}_n(C,X,d) \to V$ defined by $[u: C \to X] \mapsto u^* \Psi$.
Hence
$$ \int_{ \ol{\M}_n(C,X,d) \qu G } \ev^* {h} \cup \Eul( V) =
\int_{\ol{\M}_n(C,Z,d) \qu G} \ev^* r_{Z,G} {h} .$$
In the convex case, for sufficiently large linearization the moduli
space is empty; hence the result follows from the equality of the
wall-crossing terms \eqref{walleq} and the formula \eqref{wc1},
\eqref{wc2}.  The proof for the localized graph potentials is similar.
\end{proof} 

We apply Theorem \ref{qlef} to give a formula for a twisted localized
graph potential of moduli spaces of framed sheaves on the projective
plane.  Note that these results were already announced by
Ciocan-Fontanine-Diaconescu-Kim-Maulik, see \cite{ciocan:hilb}.  The
quantum cohomology of these moduli spaces is also the subject of work
by Maulik-Okounkov \cite{mo:qgqc}.  We recall the construction of the
moduli space of framed sheaves from Nakajima's lectures \cite[Chapter
3]{nak:lec}.  Recall the notation from \eqref{adhmpresent}. The data
$(B_-,B_+,i_-,i_+)$ forms a representation of the ADHM quiver.  The
group $S = (\C^\times)^2$ acts equivariantly on
$$ X := \End(\C^k)^{\oplus 2} \oplus \Hom(\C^r,\C^k) \oplus
\Hom(\C^k,\C^r) $$
by 
$$ (s_-,s_+)(B_-,B_+,i_-,i_+) = (s_- B_-, s_+ B_+, i_-, s_- s_+ i_+)
.$$
This preserves the locus
$$ Z_\lambda = \left\{ (B_-,B_+,i_-,i_+) \left| [B_-,B_+] + i_- i_+ =
\lambda \on{Id} \right. \right\} \subset X $$
and induces an $S$-action on the quotient $\M$.  In the case $k = 1$,
$\M$ is a deformation of the Hilbert scheme $\Hilb_k(\C^2)$ and the
$S$-action is the one induced from the $S$-action on $\C^2$.  For any
character $\chi \in \Hom(G,\C^\times) \cong \Z$, let $ \ti{\M}$ denote
Nakajima's desingularization \cite{nak:lec} of $\M$ given by
$$ \ti{\M} = Z_\lambda \qu_\chi G \subset X \qu_\chi G $$
where $\qu_\chi$ denotes the $\chi$-shifted geometric invariant theory
quotient.  Let $T \subset G$ denote the diagonal maximal torus.  We
consider the twisted Gromov-Witten theory of $X \qu_\chi G$
corresponding to the relation defining $Z_\lambda$, that is, twisted
by the Euler class of the index bundle of $\End(\C^k)$.  Although $X
\qu_\chi G$ is non-compact, the fixed points of the $S$-action are
compact:

\begin{lemma} 
\label{Sproper}
{\rm (Properness of $S$-fixed loci)}   
\begin{enumerate} 
\item The $S$-action on the moduli of parametrized stable maps
 $\ol{\M}_n(C,X \qu_\chi G,d)$ has proper $S$-fixed loci
 for any homology class $d \in H_2^{G}(X)$.
\item 
The $S$-action on the moduli space of gauged maps 
$\ol{\M}_n^{G}(C,X,d)$ has proper $S$-fixed loci 
for any homology class $d \in H_2^{G}(X)$.
\item 
The $S$-action on the moduli space of scaled gauged maps
$\ol{\M}_{n,1}^{G}(C,X,d)$ from \cite{qk3} has proper $S$-fixed loci
for any homology class $d \in H_2^{GL_n}(X)$.
\end{enumerate} 
\end{lemma} %

\begin{proof} (a)  
By the theory of symplectic resolutions discussed in
\cite{ginzburg:nakajima}, any stable map defines a parametrized
stable map to the corresponding affine quotient $X \qu G$ by
composition with the proper morphism $X \qu_\chi G \to
X \qu G$. The latter is affine with compact $S$-fixed loci,
hence any $S$-fixed stable map in $X \qu_\chi G$ projects to
an $S$-fixed point in $X \qu G$.  Since the inverse image is
proper, the claim follows.  (b) is Proposition 3.5 in
Diaconescu \cite{diac:adhm}.  However for the purposes of proving (c) we give a
different proof.  Let $(P,u)$ be an object of the fixed point substack
$\ol{\M}^{G}(C,X,d)^S$.  Thus $P \to C$ is a $G$-bundle, $u:
C \to P \times_G X$ is a section, and there exists
a homomorphism 
$$ \varphi : S \to \Aut(P) $$
such that $su = \varphi(s)u$.  After trivializing $\varphi$ at a base
point $\varphi$ defines a homomorphism still denoted $\varphi$ from
$S$ to $G$.  Let $G_\varphi$ denote the centralizer of $\varphi$.
Then $P$ admits a reduction of structure group $P_\varphi \subset P$
to $G_\varphi \subset G$.  Each $\varphi(s)$ defines an automorphism
of the associated fiber bundle $P(X)$, so that $u$ takes values in the
fixed point locus $P(X)^{\varphi} = P(X^\varphi)$.  The fixed point
locus of $S$ on $X$ defined by the homomorphism $\varphi$ is
$$ X^\varphi = \left\{ (B_-,B_+,i_-,i_+) \ \left| \ \begin{array}{l}
\Ad(\varphi(s_-,s_+)) B_\pm = s_\pm B_\pm, \\ \varphi(s_-,s_+) i_- =
i_-, \\ i_+ \varphi(s_-,s_+)^{-1} = s_+ s_- i_+ \end{array}
\right. \right\} .$$
Under the action of $S \times \C^\times$ (where $\C^\times \subset G$
is the subgroup of diagonal matrices) the subspace $X^\varphi$ splits
into a sum of subspaces with weights $(1,0,0), (0,1,0), (0,0,1),
(1,1,-1)$.  This shows the existence of a central abelian
three-parameter subgroup whose action on $X^\varphi$ has weights
contained in an open half-space.  It follows that $X^\varphi$ is
convex (take for example the one-parameter subgroup of $S \times
\C^\times$ generated by $(1,1,1)$).  So
$\ol{\M}_n^{G_\varphi}(C,X^\varphi,d)$ is compact for any class $d \in
H_2^{G_\varphi}(X)$.  Since any fixed point component arises in this
way, $\ol{\M}_n^G(C,X,d)^S$ is compact.  The argument for (c) is
similar, using that any continuous family of $S$-fixed vortices with
varying vortex parameter takes values in $X^\varphi$ for any
homomorphism $\varphi$.
\end{proof}  

By the properness results in Lemma \ref{Sproper}, the abelianization
argument for gauged potentials defined via localization at fixed point
loci of the $S$-action implies
\begin{equation} \label{walleq} \tau_{Z^\varphi,G_\varphi,\lambda,t} = \tau^V_{X,G_\varphi,\lambda,t}, \quad \forall t \in
(-1,1) .\end{equation}
Hence
$$ \tau_{Z \qu G,-} \circ \kappa_{Z,G} = \pi_T^G \tau_{X,T,-}^{V \qu
  T,\g/\t} \circ \kappa_{X,T}^{V,\g/\t} .$$
This formula can be made explicit as follows.  For any $n$-tuple of
non-negative integers $\ul{d} = (d_1,\ldots,d_k)$, $\theta = \sum c_i
\theta_i$ with $c_i \in \Z$, define
$$ \Delta_{\ul{d}}(\theta,w) := 
\frac{ \prod_{l = -\infty}^{\theta \cdot \ul{d}}  (\theta + w + l \zeta)}
{ \prod_{l = -\infty}^{0}  (\theta + w + l \zeta)} .$$
The twisted localized gauged potential for the $T$ action on
$X$ has restriction to $QH^{\leq 2}_{T}(X)$ given by
(cf. \cite{ciocan:hilb})
\begin{equation} \label{restric}
 \pi_T^G \tau_{X,T,-}^{V,\g/\t} = \sum_{d \ge 1} q^d \sum_{\ul{d}: d_1
   + \ldots + d_k = d} e^{t_0 + (t_1 (\theta_1 + d_1 \zeta) + \ldots +
   t_k (\theta_k + d_k \zeta)) / \zeta}
 \tau_{X,T,\pm}(\ul{d}) \end{equation}
where 
$$ \tau_{X,T,\pm}^{V,\g/\t}(\ul{d}) = \prod_{i \neq j}
\frac{\Delta_{\ul{d}}(\theta_i - \theta_j, \xi_- + \xi_+)
  \Delta_{\ul{d}}(\theta_i -\theta_j, 0 )} {\Delta_{\ul{d}}(\theta_i -
  \theta_j,\xi_-) \Delta_{\ul{d}}(\theta_i - \theta_j, \xi_+)}
\prod_{i =1 }^k \frac{1}{\Delta_{\ul{d}}(\theta_i,0)^r
  \Delta_{\ul{d}}(-\theta_i, \xi_- + \xi_+)^r }$$
and $\xi_-,\xi_+$ are the equivariant parameters for $S$, that is,
$H_S(\pt) = \Q[\xi_-,\xi_+]$.  The quantity $\tau_{X,T,\pm}^{(1)}$ is
studied in Konvalinka-Ciocan-Fontanine-Pak \cite{ciocan:hilb}.  The
product
$$ \prod_{i =1 }^k \frac{1}{\Delta_{\ul{d}}(\theta_i,0)^r
  \Delta_{\ul{d}}(-\theta_i, \xi_- + \xi_+)^r }$$
has leading order term $\zeta^{-1}$ only if a single $d_i $ is
non-zero and $r = 1$.  Furthermore, in the case $r =1$ the constant
term in the first product is
$$ \prod_{i \neq j} \frac{\Delta_{\ul{d}}(\theta_i - \theta_j, \xi_- +
  \xi_+) \Delta_{\ul{d}}(\theta_i -\theta_j, 0 )}
{\Delta_{\ul{d}}(\theta_i - \theta_j,\xi_-) \Delta_{\ul{d}}(\theta_i -
  \theta_j, \xi_+)} = \prod_{j \neq i} \frac{(\theta_i -\theta_j -
  \xi_-)(\theta_i - \theta_j - \xi_+)} {(\theta_i - \theta_j)(\theta_i
  - \theta_j - (\xi_- + \xi_+))} $$
It follows that for $r = 1$
\begin{multline} \tau_{X,T,\pm}^{(1)}|_{t = 0} = \sum_{d \ge 0} \frac{(-1)^d}{d} q^d
\sum_{i=1}^k (\theta_i - (\xi_- + \xi_+)) \prod_{j \neq i}
\frac{(\theta_i -\theta_j - \xi_-)(\theta_i - \theta_j - \xi_+)}
     {(\theta_i - \theta_j)(\theta_i - \theta_j - (\xi_- + \xi_+))}
     \\ \ln(1 + q) \sum_{i=1}^k (\theta_i - (\xi_- + \xi_+)) \prod_{j
       \neq i} \frac{(\theta_i -\theta_j - \xi_-)(\theta_i - \theta_j
       - \xi_+)} {(\theta_i - \theta_j)(\theta_i - \theta_j - (\xi_- +
       \xi_+))} .\end{multline}
In Konvalinka-Ciocan-Fontanine-Pak \cite{ciocan:hilb} this quantity is
equated with
$$ \exp(- \tau_{X,T,\pm}^{V,\g/\t,(1)} |_{t =0}) = (1 + q)^{k(\xi_- +
  \xi_+)}
$$
via combinatorics of Young diagrams.  It follows by \eqref{same} that
$$ \exp(\kappa_{\hat{Z},G} | _{t = 0}/\zeta) = \exp(
\tau_{X,T,\pm}^{V,\g/\t,(1)} |_{t= 0}/\zeta)
.$$
For the classes $h = t_1 h_1 + \ldots t_n h_n \in H^2_T(X)$ the value
if the potential is determined by the divisor equation, with the
result
$$ \kappa_{\hat{Z},G} = \kappa_{\hat{Z},G} | _{t = 0}(q \mapsto qe^t)
.$$
Since the classes involved in $\tau_{X,T,\pm}^{V,\g/\t}$ pair
trivially with the degrees of gauged maps, the exponential in
\eqref{restric} distributes out of the sum so that
$$ \tau_{Z \qu G,\pm} = \mu_T^G \exp( - \tau_{X,T,\pm}^{(1)}/\zeta)
\tau_{X,T,\pm} .$$
That is, letting $\delta_{ r-1}$ denote the Kronecker delta function
at $r = 1$, the twisted localized graph potential of the smoothed
moduli space of framed sheaves $Z \qu G$ restricted to $QH^{\leq 2}(X
\qu G) \cong QH^{\leq 2}_G(X)$ is
$$ \tau_{Z \qu G,-} \circ r_{Z,G} = \mu_{Z \qu G}^{X \qu T} (1 + qe^t)^{k
  \delta_{r-1}(\xi_- + \xi_+)/\zeta} \tau^{V,\g/\t}_{X,T,-}.$$
This completes the proof of Theorem \ref{desing}.

\def\cprime{$'$} \def\cprime{$'$} \def\cprime{$'$} \def\cprime{$'$}
\def\cprime{$'$} \def\cprime{$'$}
\def\polhk#1{\setbox0=\hbox{#1}{\ooalign{\hidewidtht
      \lower1.5ex\hbox{`}\hidewidth\crcr\unhbox0}}} \def\cprime{$'$}
\def\cprime{$'$}

\end{document}